\def \version {2026 -- 02 -- 10}

\documentclass[12pt]{article}

\usepackage{amsmath,amssymb,amsthm}
\usepackage[mathcal]{eucal}
\usepackage{enumerate}
\usepackage{lineno}
\usepackage{comment}

\def \bvt {\begin{verbatim}\footnotesize }
\def \evt {\end{verbatim}\footnotesize }
\newcommand{\nev}[1]{{\bf\itshape (#1.)}}

\newtheorem{thm}{Theorem}
\def \btm {\begin{thm}}
\def \etm {\end{thm}}
\newtheorem{prp}[thm]{Proposition}
\def \bpn {\begin{prp}}
\def \epn {\end{prp}}
\newtheorem{lem}[thm]{Lemma}
\def \blm {\begin{lem}}
\def \elm {\end{lem}}
\newtheorem{obs}[thm]{Observation}
\def \bob {\begin{obs}}
\def \eob {\end{obs}}
\newtheorem{rmk}[thm]{Remark}
\def \brm {\begin{rmk}}
\def \erm {\end{rmk}}
\newtheorem{cor}[thm]{Corollary}
\def \bcr {\begin{cor}}
\def \ecr {\end{cor}}
\newtheorem{con}[thm]{Conjecture}
\def \bcj {\begin{con}}
\def \ecj {\end{con}}
\newtheorem{prm}[thm]{Problem}
\def \bpm {\begin{prm}}
\def \epm {\end{prm}}
\newtheorem{dfn}[thm]{Definition}
\def \bdf {\begin{dfn}}
\def \edf {\end{dfn}}
\newtheorem{exa}[thm]{Example}
\def \bex {\begin{exa}}
\def \eex {\end{exa}}

\def \bpf {\begin{proof}}
\def \epf {\end{proof}}

\def \bsk {\bigskip}
\def \msk {\medskip}
\def \ssk {\smallskip}
\def \nin {\noindent}

\def \cA {\mathcal{A}}
\def \cB {\mathcal{B}}
\def \cC {\mathcal{C}}
\def \cE {\mathcal{E}}
\def \cF {\mathcal{F}}
\def \cG {\mathcal{G}}
\def \cI {\mathcal{I}}

\def \cX {\mathcal{X}}
\def \pp {{\mathfrak{p}}}
\def \nnn {\mathbb{N}}
\def \ovG {\overline{G}}
\def \Mad {\mathrm{Mad}}

\def \ed {\mathrm{ed}}

\def \fpr {\cF_{p,\,r}}
\def \gpr {G_{p,\,r}}

\def \gpir {G_{p_i,\,r_i}}
\def \gpjr {G_{p_j,\,r_j}}
\def \gpkr {G_{p_k,\,r_k}}
\def \summad {\Mad(G_1)+\ldots+\Mad(G_k)}
\def \summadh {\Mad(H_1)+\ldots+\Mad(H_k)}
\def \Summad {\sum_{i=1}^k \Mad(G_i)}
\def \supt {{}^{\mathsf{+}}}
\def \sL {\mathsf{L}}
\def \ssl {^{\sL}}
\def \nfty {n\to\infty}

\def \col {\mathrm{col}}
\def \ch {\mathrm{ch}}
\def \dg {\delta\supt}
\def \kaps {\kappa\supt}
\def \lams {\lambda\supt}
\def \vp {\varphi}

\newcommand{\floor}[1]{\left\lfloor #1 \right\rfloor}
\newcommand{\ceil}[1]{\left\lceil #1 \right\rceil}

\def \smin {\diagdown}
\def \es {\varnothing}

\begin{document}


\title{Nordhaus--Gaddum-type theorems for maximum average degree}
\author{Yair Caro\,\thanks{~Department of Mathematics, University of Haifa-Oranim, Tivon 36006, Israel\\
\hspace*{2em} {\tt e-mail: yacaro@kvgeva.org.il}}
 \and Zsolt Tuza\,\thanks{~HUN-REN Alfr\'ed R\'enyi Institute of Mathematics, Budapest, Hungary} $^,$\thanks{~Department of Computer Science and Systems Technology, University
of Pannonia, Veszpr\'em, Hungary\\
\hspace*{2em} {\tt e-mail: tuza.zsolt@mik.uni-pannon.hu}}}
\date{\small Latest update on \version}
\maketitle

\begin{abstract}

A $k$-decomposition $(G_1,\dots,G_k)$ of a graph $G$ is a partition of its edge set into $k$ spanning subgraphs $G_1,\dots,G_k$.
The classical theorem of Nordhaus and Gaddum bounds $\chi(G_1) +
 \chi(G_2)$ and $\chi(G_1) \chi(G_2)$ over all 2-decompositions of $K_n$.

For a graph parameter $\pp$, let $\pp(k,G) = \max \{ \sum_{i=1}^k
 \pp(G_i) \}$, taken over all $k$-decompositions of graph $G$.  
In this paper we consider $M(k,K_n) = M(k,n) = \max \{ \sum_{i=1}^k
 \Mad(G_i) \}$, taken over all $k$-decom\-positions of the complete graph
  $K_n$, where $\Mad(G)$ denotes the maximum average degree of $G$,
   $\Mad(G) = \max \{ 2e(H)/|H| : H \subseteq G \} = \max
   \{d(H) : H \subseteq G \}$.

Among the many results obtained in this paper we mention the following selected ones.

(1)\quad
 $M(k, n) < \sqrt{k}\,n$, and
  $\lim_{k\to\infty} \bigl( \liminf_{n\to\infty} \frac{M(k,n)}{\sqrt{k}\,n} \bigr)
  = 1$.

(2)\quad Exact determination of $M(2,n)$.

(3)\quad Exact determination of $M(k,n)$ when $k = \binom{n}{2} - t$,
 $0 \leq t\leq (n-1)^2\!/3$.

A main tool used along this paper is the following list variant of
 Nordhaus--Gaddum, which is of independent interest and can be adopted
 to attack other Nordhaus--Gaddum type problems, too.
A~list~$\sL$ of $k$ graphs is a multiset,
 $\sL = \{G_1,\dots,G_k\}$, where repetitions are allowed.
For a positive integer $N$ define
 $M\ssl(k,N) = \max \{ \sum \Mad(G_j) \}$ taken
  over all
   lists $\sL$ of $k$ graphs such that $\sum e(G_j) = N$.
Clearly $M(k,n) \leq M\ssl(k,\binom{n}{2})$.
We determine $M\ssl(k,N)$ for all values of the parameters involved.

(4)\quad Let $N \geq k \geq 2$ be any integers, and let $p, q, r$ be
 non-negative
 integers such that $N = k\binom{p}{2} + qp + r$, where $0 \leq q < k$ and
 $0 \leq r < p$.
Then $M\ssl(k,N) = kp - k + q$ if $r \leq (p - 1)/2$, and 
$M\ssl(k,N) = kp - k + q + 1 - 2(p - r)/(p + 1)$ if $r \geq (p-1)/2$.

Using constructions based on the theory of Steiner systems and block
 designs, we supply a plethora of cases where this list upper bound
 is indeed realized by decompositions of $K_n$ and equality
   $M(k,n)$ = $M\ssl(k,\binom{n}{2})$ holds.

Applications of these bounds to other parameters considered before
 in the literature are given.

\bsk

\nin
\hfill\break
\hspace*{2.1cm}
 Nordhaus--Gaddum theorem

\bsk

\nin

\end{abstract}

\newpage

\tableofcontents

\newpage

\section{Introduction}

Nordhaus and Gaddum \cite{NG} proved that if $G$ is a graph of 
order $n$ then  $2\sqrt n \leq \chi (G)+\chi (\ovG) \leq  n+1$ and $n \leq  
\chi (G)\cdot \chi (\ovG)\leq  (n+1)^2\!/4$.
This theorem was the starting point of a reach theory of extremal results 
of the following form, called Nordhaus--Gaddum-type theorems.

A $k$-decomposition
 of a graph $G$ is a partition of its edge set 
with $k$ spanning subgraphs
 $G_1,\dots, G_k$.
For a graph parameter $\pp$, this track of research investigates:

\msk

\begin{tabular}{cl}
$(a)$\quad & 
 $\max \, \{ \sum_{i=1}^k \pp(G_i)\}$, \\
$(b)$\quad & 
 $\min \, \{ \sum_{i=1}^k \pp(G_i)\}$, \\
$(c)$\quad & 
 $\max \, \{ \prod_{i=1}^k \pp(G_i)\}$, \\
$(d)$\quad & 
 $\min \, \{ \prod_{i=1}^k \pp(G_i)\}$,
\end{tabular}

\msk

\nin
where the maxima and minima are taken over all $k$-decompositions 
of $G$.
Analogous functions have been studied in other combinatorial structures, too.

Nordhaus and Gaddum considered the case $k = 2$ for the parameter $\pp(G) 
 = \chi (G)$, the chromatic number.
Many other graphical parameters were studied in this context since then
 (more than 2000 articles in Google Scholar).
We refer the reader to the survey \cite{AH-13}.

Also, there is an evolving literature where Nordhaus-Gaddum type theorems 
are proved for $k \geq 3$; see, for example, \cite{AW-12},
 \cite{BCP}--\cite{Bic},
 \cite{FKSSW}, \cite{HJS-12}, \cite{P-78}.

In this paper we concentrate on a problem of type $(a)$, and consider the
 extremal function
 $$
    M(k,n) := \max \, \{ \sum_{i=1}^k \Mad(G_i) \}
 $$
  taken over all $k$-decompositions of the  complete graph $K_n$,
where $\Mad(G) = \max \{ 2e(H)/|H| : H \subseteq G \}
   = \max 
\{d(H) : H \subseteq G \}$, called the maximum average degree of $G$.

Somewhat surprisingly, we were unable to find any paper concerning 
Nordhaus--Gaddum-type results for the parameter $\Mad(G)$, despite
 that it is
 used frequently in determining the
 arboricity of graphs and also in many coloring problems.
See, e.g., \cite{BK-99}, \cite{NS-22}.

The value $\Mad(G) +1$ dominates or is closely related to many  other parameters,
 including the clique number $\omega (G)$, the
 choice number or list~chromatic number
 $\ch(G)$, and the Szekeres--Wilf number---also called
  coloring~num\-ber---denoted $\col(G)$,
  which is just one greater than the
  degeneracy of $G$.
Hence $M(k,n)$ serves as an upper bound for the corresponding
 Nordhaus--Gaddum-type theorems for these parameters.
We shall use this fact to 
comment on earlier results obtained in \cite{FKSSW}.

To deal with $M(k,n)$, we also introduce a list variant of
 Nordhaus--Gaddum as follows.
A list $\sL$ of k graphs is a multiset,
 $\sL = \{G_1,\dots,G_k\}$, where repetitions are allowed.
For a positive integer $N$ define:
 $M\ssl(k,N) = \max \{ \sum \Mad(G_j) \}$ taken over all
   lists $\sL$ of $k$ graphs such that $\sum e(G_j) = N$.

Clearly $M(k,n) \leq M\ssl(k,\binom{n}{2})$.
Hence if there is a list of $k$ graphs realizing $M\ssl(k,\binom{n}{2})$
 and this list can also be realized as a decomposition of $K_n$, then
 the the inequality $M(k,n) \leq M\ssl(k,\binom{n}{2})$ turns to equality.
This idea can be adopted for other Nordhaus--Gaddum-type results
 with distinct parameters.

\subsection*{Organization of the paper}

In Section \ref{s:simple} we present simple properties of $\Mad(G)$,
 which are used later along the paper.
In Section \ref{s:general}, among several results, we prove the
 general upper bound: $M(k,n) < \sqrt{k}\,n$
 (Lemma \ref{t:conv} and Corollary~\ref{c:M-sqrt-k}).
In Section~\ref{s:k=2} we deal with cases of small $k$;
 most importantly, we determine $M(2,n)$ exactly for all $n$ (Theorem \ref{t:k=2}).
The lower-bound constructions for $M(k,n)$ for $3 \leq k \leq 7$,
 collected in Proposition \ref{p:3--7}, are closely related to those
  due to Bickle~\cite{Bic-12}.

In Section \ref{ss:main} we provide a complete determination of
 $M\ssl(k,N)$ (Theorem~\ref{t:list-k-n})
 and a sufficient condition for $M(k,n) = M\ssl(k,\binom{n}{2})$
   (Theorem \ref{t:Mad-k}), which are among the major results of this paper.
In the appendix
 we supply also a second, more algorithmic and direct proof
  of Theorem \ref{t:Mad-k}.
As an application of these results, in Section \ref{ss:highest}
 we derive the exact formula for $M(\binom{n}{2} - t,n)$
 in the entire range of $0 \leq t \leq (n-1)^2 \!/3$,
  and characterize the extremal decompositions.

In Section \ref{s:designs}
 we present various constructions based on Steiner systems,
 finite geometries, Wilson's decomposition theorem and
 pairwise balanced $K$-designs where $K$ is of the
 form $\{ p ,p+1 \}$, to obtain
 sharp or nearly sharp results on various values of $M(k,n)$.
Our methods can be applied to derive the inequality
 $(1-o_k(1))\sqrt{k}\,n - c_k < M(k, n) < \sqrt{k}\,n$ for every
 fixed $k$ and suitable constant $c_k$ (Theorem \ref{t:sqrt-k}),
 concerning the asymptotic growth of $M(k,n)$.
This is another major result of this paper.
As $k$ grows, it means that the convergence
  $\lim_{k\to\infty} \bigl( \liminf_{n\to\infty} \frac{M(k,n)}{\sqrt{k}\,n} \bigr)
  = 1$ is valid.

Lastly in Section \ref{s:appli} we consider applications of our results to several 
 parameters dominated by or related to $M(k,n)$, which were already discussed in \cite{FKSSW}, and also several further ones.
Our main result here, Theorem \ref{t:pp-k-n}, is that if
 $M(k,n) = M\ssl(k,\binom{n}{2})$, then the analogously defined extremal
 function $\pp(k,n)$ for many further graph parameters $\pp$ is equal
 to $M(k,n)$ or $\floor{M(k,n)} + k$.
This includes clique number, chromatic number, vertex/edge connectivity,
 and more.
 
Section \ref{s:concl} contains open problems for future research.

\subsection{Notation}

We apply the following standard notation for graphs $G$\,:
\,$|G|=|V(G)|$ is the number of vertices (the order of $G$),
$e(G)$ is the number of edges,
$\delta(G)$ is the minimum degree of vertices, and
$d(G) = ( \sum_{v\in V(G)} \deg(v) ) / |G| = {2e(G)}/{|G|}$
  is the average degree.

Further common notation is $\ovG$ for the complementary graph of $G$,
$G[X]$ for the subgraph induced by a vertex subset $X\subseteq V(G)$, and
$\deg(v)$ for the degree of a vertex $v\in V(G)$.

The ``support'' of $G$, denoted by $G\supt$, is defined as the subgraph
  obtained from $G$ by omitting the isolated vertices.
We write $n\supt = |G\supt|$ for the number of non-isolated verices,
 and define the ``essential average degree'' as
 $\ed(G) = ( \sum_{v\in V(G\supt)} \deg(v) ) / n\supt = 2e(G)/n\supt
  = d(G\supt)$, that is the average degree in $G\supt$.

Well known particular types of graphs are the complete graph $K_n$ of
 order $n$, and the path $P_s$ on $s$ vertices (having $s-1$ edges).
For integers $t_1,t_2>0$, the graph $t_1G_1\cup t_2G_2$ is the
 vertex-disjoint union of $t_1$ copies of $G_1$ and $t_2$ copies of $G_2$.
We also consider multisets of graphs.
The collection of $t_1$ copies of $G_1$ and $t_2$ copies of $G_2$
 will be denoted by $\{t_1*G_1,t_2*G_2\}$; and more generally,
  for multisets built upon $s$ graphs $G_1,\dots,G_s$,
  we write $\{t_1*G_1,\dots,t_s*G_s\}$.

\newpage

\subsection{Functions on spanning vs.\ non-spanning subgraphs}

Recall from above that the support $G\supt$ of a graph $G$ is defined as
 the subgraph induced by the vertices of positive degrees in $G$.  
If $G$ is the empty graph, we define $G\supt = K_1$.

\msk

For instance, if $G = K_4 \cup 4K_1$, then
 $G\supt = K_4$, $d(G) = 12/8 = 1.5$, and $d(G\supt) = 3 = \ed(G)$.
As another parameter, considering edge connectivity $\lambda(G)$,
 for the same $G = K_4 \cup 4K_1$ we have
 $\lambda(G) = 0$ while $\lambda(G\supt) = \lambda(K_4) = 3$.  

\msk

Let $H$ be any given host graph, and $\pp$ any given graph invariant.
Let us introduce the following functions.

\msk

(1)\quad $\pp(k,H) = \max \{ \sum \pp(G_j) \}$ taken over all
 decompositions of $E(H)$ into $k$ spanning subgraphs of $H$.
Note that isolated vertices are allowed in any of the $G_j$.
If $H = K_n$, we write $\pp(k,n)$.

\msk

(2)\quad $\pp\supt(k,H) = \max \{ \sum \pp(G_j\supt) \}$ taken over all
 decompositions of $E(H)$ into $k$ spanning subgraphs of $H$.
Note that $G\supt$ eliminates the isolated vertices from $G$.
If $H = K_n$, we write $\pp\supt(k,n)$.

\msk

A list $\sL$ of k graphs is a multiset,
 $\sL = \{G_1,\dots,G_k\}$, where repetitions are allowed.

\msk

(3)\quad For a positive integer $N$ define:
 $\pp\ssl(k,N) = \max \{ \sum \pp(G_j) \}$ taken
  over all\footnote{~Interesting problems arise also where the universe
   of graphs is restricted, or just the admissible lists themselves are
   required to satisfy certain constraints.
    However, in this paper we deal with general lists
   admitting all graphs.}
   lists $\sL$ of $k$ graphs such that $\sum e(G_j) = N$.
It means that $N$ edges with non-specified endpoints are given,
 and we are asked to build $k$ graphs $G_1,\dots,G_k$ from them
 to maximize $\pp(G_1)+\ldots+\pp(G_k)$.

\msk

It depends on the invariant $\pp$ whether or not it is sensitive
 with respect to the presence of isolated vertices.
In the current examples $d(G)$ and $\lambda(G)$ are, but $\Mad(G)$
 and $\ed(G)$ aren't.
For insensitive parameters it is not necessary to assume that all
 the graphs $G_1,\dots,G_k$ forming a decomposition of $K_n$ are
 spanning subgraphs, the definitions can be relaxed to subgraphs.
Further, if $\pp$ is monotone increasing in the sense that
 $\pp(G)\leq \pp(H)$ whenever $G\subset H$, then edge-disjoint packings
 may be studied instead of decompositions.

\section{Some simple properties of Mad({\itshape G})}
\label{s:simple}

We introduce the following term, which will simplify some arguments.

\bdf   \nev{Free edge}
In a graph\/ $G$, we call an\/ $e\in E(G)$ a free edge if
 its removal does not modify the\/ $\Mad$ value; i.e.,\/
 $\Mad(G-e) = \Mad(G)$.
\edf

\bpn   \nev{Monotonicity}
\label{p:monot}
For every\/ $n$ and\/ $k$,
 \begin{itemize}
  \item[$(i)$] $M(k,n+1) \geq M(k,n) + 1$ if\/ $k\leq \binom{n}{2}$.
  \item[$(ii)$] $M(k+1,n) \geq M(k, n) + 1/3$ if\/ $k < \binom{n}{2}$.
 \end{itemize}
\epn

\bpf
Let $G_1,\dots,G_k$ decompose $K_n$ with $M(k,n) = \Summad$, and let
 $V(K_{n+1}) = V(K_n) \cup \{w\}$.
Assume that $\Mad(G_k)$ is attained on a vertex set $X$, i.e.,
 $\Mad(G_k) = 2e(G_k[X])/|X|$.
Extend now $G_k$ with vertex $w$ and all edges $vw$ where $v\in V(K_n)$,
 hence the join $G_k + K_1$.
In the extended $G_k$ the number of edges induced by $X \cup \{w\}$ is
 $e(G_k[X])+|X|$, therefore the maximum average degree is increased to
 $$
   \frac{2e(G_k[X])+2|X|}{|X|+1} = \frac{2e(G_k[X])+(|X|-1)}{|X|+1} +1
    \geq \frac{2e(G_k[X])}{|X|} +1
 $$
  because $X$ cannot induce more than $\binom{|X|}{2}$ edges.
This proves $(i)$.

For $(ii)$, let
 $\Summad = M(k,n)$, and suppose that $e(G_1)$ is largest among the $e(G_i)$.
Then $e(G_1)>1$, hence the number $n_1$ of non-isolated vertices
 in $G_1$ is at least 3.
If $G_1$ has a free edge $e$, then taking $G_{k+1}=K_2$ whose edge is $e$
 we increase the sum of the $\Mad(G_i)$ by 1.
If this is not the case, then $\Mad(G_1) = 2e(G_1)/n_1$, each edge
 contributing to $\Mad(G_1)$ with $2/n_1\leq 2/3$.
Thus, picking an arbitrary edge $e$ and creating also here $G_{k+1}=K_2$
 whose edge is $e$, $\Mad(G_1) - \Mad(G_1-e)\leq 2/3$ holds, while
 $\Mad(G_{k+1})=1$.
So, this modification increases the sum of the $\Mad$ values by at least 1/3.
\epf

Comparison of the two parts of Theorem \ref{t:PSTS} provides us with
 $(n-1)^2\!/6$ values of $k$ for which the increase of 1/3 given in
 $(ii)$ above is tight.

\bpn
\label{p:conn-ind}
If\/ $G$ is not empty, then there is an induced connected
 subgraph\/ $H\subseteq G$ realizing\/ $\Mad(G)=d(H)$.
\epn

\bpf
Assume that $H\subseteq G$ satisfies $d(H)=\Mad(G)$, and under
 this condition $|V(H)|$ is smallest.
It is clear that all edges induced by $V(H)$ belong to $E(H)$.
Suppose for a contradiction that $H$ is disconnected, say $H=H'\cup H''$
 and there is no edge between $H'$ and $H''$.
We have $d(H')=2e(H')/|H'|$, $d(H'')=2e(H'')/|H''|$, and
 $d(H)=2(e(H')+e(H''))/(|H'|+|H''|)$.
Assuming $d(H') \leq d(H'')$, and taking their weighted average
 we clearly have
 $$
   d(H') \leq \frac{|H'|}{|H'|+|H''|} \, d(H')
    + \frac{|H''|}{|H'|+|H''|} \, d(H'') \leq d(H'') \,.
 $$
Observe that the expression in the middle is exactly $\Mad(H)$.
Thus, $H''$ has fewer vertices than $H$ has, while its
 average degree is not smaller.
This contradicts the choice of $H$.
\epf

\bpn
\label{p:tree}
If\/ $T$ is a tree, then\/ $\Mad(T) = \ed(T) = d(T) = 2 - 2/|T|$.
More generally, if\/ $G$ is a maximal\/ $\delta$-degenerate graph
 on more than\/ $\delta$ vertices, then\/
 $\Mad(G) = \ed(G) = d(G) = 2\delta - \frac{\delta(\delta+1)}{|G|}$.
\epn

\bpf
Assume that the subgraph $G'\subseteq G$ attains $d(G') = \Mad(G)$, and
 $G'$ contains no free edges.
We apply the facts that a maximal $\delta$-degenerate graph
 on $n>\delta$ vertices has exactly
 $\delta n - \binom{\delta+1}{2}$ edges, moreover every subgraph
 of a $\delta$-degenerate graph is also $\delta$-degenerate.
In this way we obtain
 \begin{eqnarray}
   \Mad(G) & = & d(G') \ = \ \frac{2e(G')}{|G'|} \ \leq \
    2\delta - \frac{\delta(\delta+1)}{|G'|} \nonumber \\
     & \leq & 2\delta - \frac{\delta(\delta+1)}{|G|} \ = \ d(G) \ = \ \ed(G) \nonumber \\
     & \leq & \Mad(G) \,. \nonumber
 \end{eqnarray}
Thus, equality holds throughout.
It also follows that $G'=G$.
\epf

Concerning lower-bound constructions for $M(k, n)$ with fixed $k$,
 it is useful to state the following simple implication explicitly.

\bpn
\label{p:split}
If\/ $p>q>1$ and\/ $G$ is the complete split graph\/ $K_p - K_q$, then\/
 $\Mad(G) = \ed(G) = d(G)$.
\epn

\bpf
The graph $K_p - K_q$ is $(p-q)$-degenerate, and it is maximal
 with this property because inserting a new edge creates a
 subgraph $K_{p-q+2}$.
Thus, Proposition \ref{p:tree} can be applied.
\epf

\bob
Since\/ $\Mad(G)\geq\ed(G)\geq d(G)=2e(G)/|V(G)|$ holds for every graph\/ $G$,
 the following inequalities are immediately seen, where the maximum is
 taken over all\/ \emph{edge-disjoint packings} $\cG = \{G_1,\dots,G_k\}$
 of\/ $k$ subgraphs of\/ $K_n$\,:
  $$
    M(k,n) = \max_{\cG} \Summad \geq \max_{\cG} \sum_{i=1}^k \ed(G_i)
     \geq \max_{\cG} \sum_{i=1}^k d(G_i) \,.
  $$
\eob

\bob
To have\/ $\Mad(G)=d(G)$ for a graph\/ $G$, also\/ $\Mad(G)=\ed(G)$ must hold.
This forces\/ $\delta(G)\geq d(G)/2$, for otherwise the removal of a vertex\/
 $v$ of smallest degree\/ $d(v)=\delta(G)<|V(G)|/2$ would yield\/
  $d(G - v) = 2(e(G) - \delta(G))/(|V(G)|-1) > (2e(G) - d(G))/(|V(G)|-1)
  = 2e(G)/|V(G)| = d(G)$.
Dense examples of graphs\/ $G$ with\/ $\Mad(G)=d(G)$ are obtained, e.g.,
 by removing at most\/ $n/2$ edges from\/ $K_n$.
Sparser examples are all regular graphs, and many more.
\eob


\section{General estimates and constructions}
\label{s:general}

In the first part of this section we determine the largest possible
 value of $\Mad(G)$ among graphs with a given number $m\geq 1$
 of edges; that is,
  $$
    g(m) := \max \{ \Mad(G) \mid e(G)=m \} \,.
  $$
As one may observe, with the notation defined in the introduction, we have
 $g(m) = M\ssl(1,m)$.

In the second part of the section we describe a structure class of edge
 decompositions, on which it would suffice to determine
 $\max\{\Summad\}$ in order to determine $M(k,n)$.

\subsection{Single graph optimizing Mad({\itshape G\/})}

\btm
\label{t:max-g-m}
Let\/ $p\geq 2$ be an integer.
 \begin{itemize}
  \item[$(i)$] If\/ $m=\binom{p}{2}$, then\/ $g(m)=p-1$, with the
   unique extremal graph\/ $K_p$.
  \item[$(ii)$] If\/ $\binom{p}{2} < m < \binom{p+1}{2}$,
    say $m=\binom{p}{2}+r$ where $0<r<p$, then\/
   $g(m)=\max\{p-1, 2m/(p+1)\}$, attained by precisely the
   following graphs\/ $G$ having\/ $m$ edges:
    \begin{enumerate}[(a)]
  \item any\/ $G$ with\/ $K_p\subset G$, if\/ $0<r<\frac{p-1}{2}$;
  \item any\/ $G$ with\/ $|V(G)|=p+1$, if\/ $\frac{p-1}{2}<r<p$;
  \item both types (a) and (b), if\/ $r=\frac{p-1}{2}$.
    \end{enumerate}
 \end{itemize}
\etm

\bpf
Let $G$ be any graph with $m$ edges, and $X\subseteq V(G)$ a set
 of vertices inducing a subgraph $H$ with $d(H)=\Mad(G)$.
If $|X|\leq p$, then $X$ can induce at most $\binom{|X|}{2}$
 edges, thus $d(H) = 2e(H)/|X| \leq |X|-1 \leq p-1$.
If $|X|\geq p+1$, then $d(H) = 2e(H)/|X| \leq 2e(G)/(p+1)$.

Both bounds are attained by graphs of $m$ edges in the given range,
 namely by $K_p$ supplemented with arbitrarily positioned
 $m-\binom{p}{2}$ further edges, and by all graphs with
 $p+1$ vertices and $m$ edges, respectively.
\epf

\bdf   \nev{Family of extremal (p,r)-graphs}
For\/ $p\geq 2$ and\/ $p\geq r\geq 0$, let\/
 $m = \binom{p}{2} + r$.
Based on the above, we denote:

{\small
 $$
   \fpr =  \begin{cases}
   \begin{tabular}{ll}
    $\{K_p\}$ & if\/  $r=0$\,; \\
    $\{ G \mid e(G) = m$\,;\, $K_p\subset G \}$
      & if\/  $0 < r < (p-1)/2$\,; \\
    $\{ G \mid e(G) = m$\,;\, $K_p\subset G$
     or $|G| = p+1 \}$
      & if\/  $r = (p-1)/2$\,; \\
    $\{ G \mid e(G) = m$\,;\, $|G| = p+1 \}$
      & if\/  $r > (p-1)/2$\,; \\
    $\{K_{p+1}\}$ & if\/  $r = p$\,.
   \end{tabular}
    \end{cases}
 $$
 }
\edf

\bdf \nev{Representative}
For\/ $r = 0$ the representative\/ $G_{p,\,0}$ of\/ $\cF_{p,\,0}$ is\/ $K_p$.
For any\/ $p\geq r\geq 1$, we denote by\/ $\gpr$ the graph with\/
 $p+1$ vertices, obtained from\/ $K_{p}$ by
 supplementing it with a further vertex of degree\/ $r$.
We view this\/ $\gpr$ as the representative of the family\/ $\fpr$ of extremal
 graphs for\/ $g(m)$ where\/ $m=\binom{p}{2}+r$.
In particular,\/ $G_{p,\,p}\cong K_{p+1}$.
\edf

As in the definition of representative we allow both $r =0$ and $r = p$,
 we can view $K_{p+1}$ as $G_{p+1,\,0}$ a representative of
 $\cF_{p+1,\,0}$ and also as $G_{p,\,p}$ a representative of $\cF_{p,\,p}$,
 where the two families coincide.

One may also note that, for $0<r<(p-1)/2$, the $r$ edges incident with the
 low-degree vertex of $G_{p,\,r}$ are free edges (and the others are not).

\newpage

\subsection{Equivalent condition for {\itshape M\hspace{1pt}}({\itshape k}\hspace{1pt},\,{\itshape n\/})}
\label{ss:equiv}

Concerning the next result, let us recall that the graphs $G_1,\dots,G_k$
 attaining $M(k,n)=\summad$ do not necessarily decompose $K_n$,
 some $e\in E(K_n)$ may not belong to any of those $G_i$.

\bpn
\label{p:X1--Xk}
Let\/ $n\geq 3$ and\/ $k\geq 2$.
Suppose that the packing of graphs\/ $G_1,\dots ,G_k$
 optimzes\/ $M(k,n)$.  
Let\/ $X_j  \subset V(G_j )$, $j = 1,\dots,k$ be a subset on which\/
 $\Mad(G_j)$ is realized, namely\/ $2e(G_j[X_j])/ |X_j|  = \Mad(G_j)$.  

Assume\/ $|X_1| \leq \dots \leq |X_k|$.
Then, taking\/ $H_1$ as the complete graph on\/ $X_1$ and\/ $V(H_j)=X_j$,
  $E(H_j) =  E(K_n[X_j]) \smin (\,\cup_{i<j} E(H_i))$ for\/ $j\geq 2 $,
 also the graphs\/ $H_1,\dots , H_k$ optimize\/ $M(k,n)$ such that\/
 $\Mad(H_j)$ is realized on\/ $X_j$.
\epn

\bpf
To simplify notation, let us write $\cE = E(G_1)\cup\cdots\cup E(G_k)$
 and $\cX = \binom{X_1}{2}\cup\cdots\cup\binom{X_k}{2}$, both viewed as
 subsets of $E(K_n)$.

If $e\in E(K_n)\smin\cX$, then $e$ plays no role in $M(k,n)$, because
 $\Mad(G_j)=\Mad(G_j[X_j])$ for every $j$.
Hence, we may assume $\cE\subseteq \cX$ without loss of generality.
On the other hand, if $e\in\cX\smin\cE$ held for an edge $e$, then there
 would exist $j$ with $e\in\binom{X_j}{2}$ and re-defining
 $E(G_j) := E(G_j)\cup\{e\}$ would increase $\Summad$ by $2/|X_j|$,
  contradicting the assumption that $G_1,\dots,G_k$ realizes $M(k,n)$.
Thus $\cX\subseteq\cE$, consequently $\cX=\cE$.

We will prove that the re-distribution of the edges belonging to $\cX$ by the rules $E(H_1) = \binom{X_1}{2}$ and
 $E(H_j) = \binom{X_j}{2} \smin (\,\cup_{i<j} \binom{X_i}{2})$ for $j\geq 2$
 yields $\summadh\geq\summad$, which clearly implies that also
 $H_1,\dots,H_k$ optimize $M(k,n)$.
For this purpose the following facts can be observed for any
 $1\leq i<j\leq k$ and any edge $e \in \binom{X_i}{2} \cap \binom{X_j}{2}$.
 \begin{itemize}
  \item If $|X_i|<|X_j|$, then $e\in E(G_j)$ is impossible, because
   in that case moving $e$ from $G_j$ to $G_i$ would increase
    $\Mad(G_i)+\Mad(G_j)$ with $2/|X_i|-2/|X_j|>0$, a contradiction to
    $\summad = M(k,n)$.
  \item If $|X_i|=|X_j|$, then for $e$ it is irrelevant whether
   $e\in E(G_i)$ or $e\in E(G_j)$, because in either case the
    contribution of $e$ to $\Mad(G_i)+\Mad(G_j)$ is $2/|X_i|=2/|X_j|$.
   (If moving $e$ from $G_j$ to $G_i$ modifies $G_j$ to a graph in which
    $X_j$ does not represent $\Mad$ anymore, it means that $\summad$
    is improvable, contradicting that $G_1,\dots,G_k$ optimize $M(k,n)$.)
 \end{itemize}
The first property implies that, already at the beginning, each $e\in\cX$
 belongs to a $G_j$ such that $|X_j| = \min\{|X_i| \mid i\in\cI_e\}$.
And then, if this minimum is attained for more than one subscript $j$,
 the second property yields that we are free to choose any such $j$,
 specifically we can prescribe $j=\min\{i\in\cI_e\}$ for $e$
 because $|X_1|\leq \cdots \leq |X_k|$.
Applying this ``shift to the left'' (or First Fit) rule,
 $j=\min\{i\in\cI_e\}$ for all $e\in\cX$ exactly means
 for the obtained graphs $H_1,\dots,H_k$ that $H_1$ is the complete graph
 on the vertex set $X_1$, and for every $j>1$ the graph $H_j$ contains
 exactly those edges from $\binom{X_j}{2}$ which do not occur in
 $H_1\cup\cdots\cup H_{j-1}$.
\epf

It is also valid that the $\Mad$-realizing sets for an
 optimal decomposition cover the entire $V(K_n)$.

\bpn
\label{c:unio-X}
Let\/ $n\geq 3$ and\/ $k\geq 2$.
Suppose that the packing of graphs\/ $G_1,\dots ,G_k$ optimzes\/ $M(k,n)$,
 and the sets\/ $X_j  \subset V(G_j )$ for $j = 1,\dots,k$ realize\/
 $\Mad(G_j)$, i.e.,\/ $2e(G_j[X_j])/ |X_j|  = \Mad(G_j)$.
Then\/ $X_1\cup\cdots\cup X_k = V(K_n)$.
\epn

\bpf
Otherwise choose any vertex $v\in V(K_n)\smin (X_1\cup\cdots\cup X_k)$.
The edges $vx\in E(K_n)$ for $x\in X_1$ have no contribution to $\summad$,
 therefore we can replace $G_1$ with the graph
  $H := ( X_1\cup\{v\}, E(G_1[X_1] \cup \{ vx \mid x\in X_1 \} )$
 to obtain a modified packing of $k$ subgraphs in $K_n$.
Then, denoting $e_1 = e(G_1[X_1])$ and $x_1 = |X_1|$, and observing that
 $e_1\leq \frac{1}{2} (x_1^2 - x_1)$, we obtain
 $$
   \Mad(H) - \Mad(G_1) = \frac{2e_1+2x_1}{x_1+1} - \frac{2e_1}{x_1}
    = \frac{2x_1^2 - 2e_1}{x_1^2 + x_1} \geq 1 \,.
 $$
This contradicts the assumption that $\summad$ attains the maximum, $M(k,n)$.
\epf

From Propositions \ref{p:X1--Xk} and \ref{c:unio-X} we also obtain:

\bcr
The determination of\/ $M(k,n)$ is equivalent to the following
 optimization problem:
 $$
     \sum_{i=1}^k \frac{2e(G_i)}{|X_i|} \qquad \longrightarrow \qquad \max
 $$
  under the following conditions:
 \begin{itemize}
  \item $X_1\cup\cdots\cup X_k = V(K_n)$;
  \item $V(G_i)=X_i$ for $i=1,\dots,k$;
  \item $E(G_j) =
    \binom{X_j}{2} \smin\, {\displaystyle \cup_{i<j}} \binom{X_i}{2}
     $.
 \end{itemize}
 \ecr

An advantage of this formulation is that it admits an integer
 programming approach, which may be efficient when $k$ and $n$ are
 not too large.

\subsection{Upper bound for {\itshape M\hspace{1pt}}({\itshape k}\hspace{1pt},\,{\itshape n}) from continuous relaxation}

\blm
\label{t:conv}
Suppose\/ $x_1,\dots,x_k > 0$   are reals such that\/
 $\sum_{j=1}^k x_j = \binom{n}{2}$, and let\/ $y_1,\dots,y_k>0$ be
 the reals with\/ $y_ j(y_ j - 1)/2 = x_ j$. 
Then:
 \begin{itemize}
  \item[$(1)$]
   $\max\, \{\,\sum_{j=1}^k (y_j - 1)\}$ is obtained when
    all\/ $x_j$ are equal, namely\/ $x_ j = \binom{n}{2}/k$
     for\/ $j=1,\dots,k$.
  \item[$(2)$]
   $\max\, \{\,\sum_{j=1}^k (y_j - 1)\} = (\sqrt{k^2 + 4kn^2 - 4kn} - k)/2
    < \sqrt{k}\, n$.
  \item[$(3)$]
   If\/ $k = \binom{n}{2}/\binom{t}{2}$, then\/ $y_j(y_j-1)/2 = t(t-1)/2$
    and the maximum equals\/ $k(t-1) = n(n-1)/t$.
 \end{itemize}
\elm

\bpf
The definition of $y_j$ yields $y_j^2 - y_j - 2x_j = 0$, the positive
 root is $y_j = \frac{1}{2} ( 1 + \sqrt{8x_j+1}\, )$.
Since $a\sqrt{x}+b$ is a concave function for $a>0$, the assertion of (1)
 follows by Jensen's inequality.

Plugging in $x=\binom{n}{2}/k$ for all $x_j$, we obtain
$y(y-1)/2 = n(n-1)/2k$, hence $ky^2 - ky - n^2 +n = 0$
and $y= (k + \sqrt{k^2 + 4kn^2 - 4kn})/2k$ for all $y_j$, therefore
 \begin{eqnarray}
   k(y-1) & = & k((k + \sqrt{k^2 + 4kn^2 - 4kn} - 2k)/2k \nonumber \\
   & = & (\sqrt{k^2 + 4kn^2 - 4kn} - k))/2 \nonumber \\
   & < & \sqrt{k}\,n \nonumber
 \end{eqnarray}
since the last two lines are positive and eliminating the first square root by squaring we get $4kn^2 +k^2 - 4kn\leq (2\sqrt{k}\,n +k) ^2 = 4kn^2 + 4k^{3/2}\,n +k^2$, which is obviously true.
This proves (2).

Finally, (3) is a consequence of (1).
Indeed, by (1) the maximum of $\sum_{j=1}^k (y_j - 1)$
 is attained when
 $
   \binom{y_j}{2} = x_j = \binom{n}{2} / k = \binom{t}{2}
 $
holds for all $j=1,\dots,k$\,; consequently $y_j = t$ and
 $$
   (y_1-1) + \dots + (y_k-1) = (t-1)\cdot k =
     (t-1) \binom{n}{2} / \binom{t}{2} = n(n-1)/t \,.
 $$

\epf

\brm
\label{r:sqrt}
Literally the same computation yields:
$$
  M\ssl(k,N) < \sqrt{2kN} \,.
$$
\erm

It is worth noting that the jump of estimate in the last step yielding
 $\sqrt{k}\,n$ in the proof of Lemma \ref{t:conv}
  is as large as $cn^2$ when $k$ is proportional to $n^2$,
 hence (2) cannot be within $(1+o(1))$ to the optimum if $k$ is
 quadratic in $n$.

In the same flavour, the following real relaxation of
 Theorem \ref{t:max-g-m} can also be stated.

\bpn
\label{p:gm-conv}
For every\/ $m>0$, $g(m)\leq y-1$, where\/ $y>0$ is the real number
 satisfying\/ $y(y -1)/2 = m$.
\epn

\bpf
Let $m=\binom{p}{2}+r$, where $p$ and $r$ are integers, and $0\leq r<p$.
We have seen that $g(m) = \max\{p-1,2m/(p+1)\}$ holds.

If $g(m)=p-1$, then $g(m)+1 = p \leq y$ because
 $\binom{p}{2} \leq m = \binom{y}{2}$, with equality if and only if $r=0$.

In the other case we have $g(m) = 2m/(p+1) = y(y-1)/(p+1) < y-1$
 with strict inequality, because $\binom{y}{2} = m < \binom{p+1}{2}$
 and therefore $y<p+1$.
\epf

\bcr
\label{c:M-sqrt-k}
$M(k,n) \leq M\ssl(k,N) < \sqrt{2kN} < \sqrt{k}\, n$, where\/
 $N=\binom{n}{2}$.
\ecr

\bpf
Let $G_1,\dots, G_k$ be an edge decomposition of $K_n$, and let
 $y_j(y_j-1)=2e(G_j)$ for $j=1,\dots,k$.
Applying Proposition \ref{p:gm-conv} and part (2) of Lemma \ref{t:conv}
 we obtain: $M(k,n) = \Summad \leq \sum_{j=1}^k (y_j-1) < \sqrt{k}\, n$.
The intermediate $\sqrt{2kN}$ comes from Remark \ref{r:sqrt}.
\epf

Later in Theorem \ref{t:sqrt-k} we shall see that the coefficient
 $\sqrt{k}$ of $n$ is very close to be best possible, for every $k$.

\newpage

\section{{\itshape M\hspace{1pt}}({\itshape k}\hspace{1pt},\,{\itshape n}) for small {\itshape k}}
\label{s:k=2}

In this section we consider the cases $2\leq k \leq 7$.
The exact value of $M(2,n)$ will be determined for all $n\geq 3$,
 and lower bounds will be given for the other values of $k$.

\subsection{Exact formula for {\itshape k\/} = 2}
\label{ss:k=2}

\btm
\label{t:k=2}
If\/ $G$ is any graph of order\/ $n\geq 3$, then
$$
  \Mad(G) + \Mad(\ovG) \leq
 \begin{cases}
  \begin{tabular}{ll}
   $(5n^2 - 6n +1)/4n$ & if\/ $n$ is odd, \\
   $(5n^2 - 6n)/4n$ & if\/ $n$ is even.
  \end{tabular}
 \end{cases}
$$
Moreover, both bounds are tight, hence the above formulas give the
 right value\/ $M(2,n)=5n/4-3/2+\epsilon_n$, with\/ $\epsilon_n=0$ if\/ $n$
 is even and\/ $\epsilon_n=\frac{1}{4n}$ if\/ $n$ is odd, for every\/
  $n\geq 3$.
\etm

\bpf
Let $X$ and $Y$ be subsets of $V(K_n)$, such that $\Mad(G)=d(G[X])$ and
 $\Mad(\ovG)=d(\ovG[Y])$.
Say, $|X|=x$ and $|Y|=y$.
We may assume $x\leq y$ without loss of generality, as
 we are free to switch between $G$ and $\ovG$ if necessary.
Also, from Corollary \ref{c:unio-X} applied for $k=2$
 (that is, $G_1=G$ and $G_2=\ovG$),
 we see that $\Mad(G) + \Mad(\ovG)$ can be largest only if
 $X\cup Y = V(K_n)$, hence $x+y\geq n$ and $|X\cap Y| = x+y-n$.
Moreover, by Proposition \ref{p:X1--Xk}, we may restrict
 attention to the case $E(G[X]) = \binom{X}{2}$ and $E(\ovG[Y]) =
 \binom{Y}{2}\smin \binom{X\cap Y}{2}$, whereas the edges (if any)
 connecting $X\smin Y$ with $Y\smin X$ can be distributed arbitrarily
 between $G$ and $\ovG$ and do not have any contribution to
 $\Mad(G) + \Mad(\ovG)$.
Consequently,

 \begin{eqnarray}
   \Mad(G) + \Mad(\ovG) & \leq &
    x-1 + \frac{y(y-1)-(x+y-n)(x+y-n-1)}{y} \nonumber \\
     & = &
      \frac{2nx + 2ny + x - x^2 - y - xy - n^2 - n}{y} \nonumber \\
     & = &
      2n - 1 - x + \frac{2nx + x - x^2 - n^2 - n}{y} \nonumber \\
     & \leq &
      2n - 1 - x + \frac{2nx + x - x^2 - n^2 - n}{n} \nonumber \\
     & = &
      n - 2 + \frac{x(n + 1 - x)}{n} \nonumber \\
     & \leq &
      n - 2 +
       \frac{1}{n} \floor{\frac{(n + 1)^2}{4}} ,
        \nonumber
 \end{eqnarray}
  which is equivalent to the stated upper bound.
Tightness is shown by the graph $G=K_{\floor{(n + 1)/2}}$ or
 $G=K_{\ceil{(n + 1)/2}}$ (which coincide if $n$ is odd)
 and its complement.
\epf

\subsection{Lower-bound constructions for small {\itshape k}}
\label{ss:small}

In the following small cases of $k$ we describe two constructions, which
 are substantially different.
Although the corresponding sums of the $\Mad(G_i)$ are not exactly the
 same, the differences are marginal.
However, there seems to be no smooth transition between the two
 structures presented for the same~$k$.
This fact may be viewed as an indication that the precise determination
 of $M(k,n)$ for these values of $k$ may be not at all easy.

\bpn
\label{p:3--7}
There exist constants\/ $c_3,\dots,c_7$ such that
 \begin{itemize}
  \item $M(3,n) \geq 3n/2 - c_3$,
  \item $M(4,n) \geq 5n/3 - c_4$,
  \item $M(5,n) \geq 11n/6 - c_5$,
  \item $M(6,n) \geq 2n - c_6$,
  \item $M(7,n) \geq 7n/3 - c_7$.
 \end{itemize}
\epn

\bpf
Let us recall from Proposition \ref{p:X1--Xk} that there is a natural
 two-way correspondence between $\Mad$-optimal edge decompositions
 $G_1,\dots,G_k$ of $K_n$ and sequences $X_1,\dots,X_k$ of subsets of
 $V(K_n)$ arranged in increasing order of cardinalities.
More explicitly, the graphs $G_j$ can be obtained by the recursive rule
 $E(G_j) = \binom{X_j}{2} \smin ( \, \cup_{i<j} \, E(G_i) \, )$.

In each case of $k=3,4,5,6$ we create two sequences $G_1,\dots,G_k$.
In one of them we begin with three sets $A_1,A_2,A_3$, mutually
 disjont, of respective cardinalities $\floor{n/3}\leq a_1\leq a_2\leq
  a_3\leq \ceil{n/3}$, with $a_1+a_2+a_3=n$.
In the other, which is usually better with a marginal extent,
 we begin with two disjont sets $B_1,B_2$
 of cardinalities $|B_1|=b_1=\floor{n/2}$ and $|B_2|=b_2=\ceil{n/2}$.

\msk

\nin
\underline{$k=3$\,:}\quad
Set $X_1 = A_1 \cup A_2$, $X_2 = A_1 \cup A_3$, $X_3 = A_2 \cup A_3$.
Then $G_1\cong K_{a_1+a_2}$, $G_2\cong K_{a_1+a_3}-K_{a_1}$,
 $G_3\cong K_{a_2,a_3}$; hence $\Mad(G_1)=a_1+a_2-1\sim 2n/3$,
 $\Mad(G_2)\sim n/2$, $\Mad(G_3)\sim n/3$, summing up to $3n/2$
 approximately.

Alternatively, set $X_1=B_1$, $X_2=B_2$, $X_3=B_1\cup B_2=V(K_n)$.
Then $G_1\cong K_{\floor{n/2}}$,
 $G_2\cong K_{\ceil{n/2}}$ (vertex-disjoint), and
  $G_3\cong K_{\floor{n/2},\ceil{n/2}}$.
Here each $G_i$ has average degree near to $n/2$.
\msk

\nin
\underline{$k=4$\,:}\quad
Let $X_i = A_i$ for $i=1,2,3$ (hence $G_i\cong K_{a_i}$), and
 $X_4 = V(K_n)$ (hence $G_4\cong K_{a_1,a_2,a_3}$).
Then $\Mad(G_1)+\Mad(G_2)+\Mad(G_3)=n-3$, and $\Mad(G_4)\sim 2n/3$.
Interestingly, the sets $A_i$ also offer the solution
 $X_1=A_1$, $X_2=A_1\cup A_2$, $X_3=A_1\cup A_3$, $X_4=A_2\cup A_3$.
Then $G_1\cong K_{a_1}$, $G_2\cong K_{a_1+a_2}-K_{a_1}$,
 $G_3\cong K_{a_1+a_3}-K_{a_1}$, $G_4\cong K_{a_2,a_3}$.
The corresponding respective average degrees are about $n/3$, $n/2$,
 $n/2$, $n/3$.

The other approach starts with $X_1=B_1$, $X_2=B_2$, i.e., $G_1\cong
 K_{b_1}$ and $G_2\cong K_{b_2}$, hence $\Mad(G_1)+\Mad(G_2)=n-2$.
Then split $X_2$ into two parts, $X_2=Y\cup Y'$, $|Y|=\floor{b_2/2}$,
 $|Y'|=\ceil{b_2/2}$; and set $X_3=X_1\cup Y$, $X_4=X_1\cup Y'$.
The obtained graphs are $G_3\cong K_{b_1,\floor{\floor{b_2/2}/2}}$ and
 $G_4\cong K_{b_1,\ceil{\floor{b_2/2}/2}}$, each having its average
 degree nearly $n/3$.

\msk

\nin
\underline{$k=5$\,:}\quad
Let $X_1=A_1$, $X_2=A_2$, $X_3=A_1\cup A_2$, $X_4=A_1\cup A_3$,
 $X_5=A_2\cup A_3$.
This yields $G_1\cong K_{a_1}$, $G_2\cong K_{a_2}$, $G_3\cong K_{a_1,a_2}$,
 $G_4\cong K_{a_1+a_3}-K_{a_1}$, $G_5\cong K_{a_2,a_3}$.
The average degrees are nearly $n/3$, $n/3$, $n/3$, $n/2$, and $n/3$,
 respectively.

Also here we can alternatively start with $X_1=B_1$, $X_2=B_2$,
 hence generating $\Mad(G_1)+\Mad(G_2)=n-2$.
For the other three graphs we cut both $B_1$ and $B_2$ into half,
 say $B_1=B_1'\cup B_1''$ and $B_2=B_2'\cup B_2''$, each of those four
 sets having cardinality $\floor{n/4}$ or $\ceil{n/4}$.
Then set $X_3=B_1'\cup B_2'$, $X_4=B_1''\cup B_2'$, $X_5=B_1\cup B_2''$.
The average degree in $G_3$ and in $G_4$ is about $n/4$, and in $G_5$
 it is about $n/3$.

\msk

\nin
\underline{$k=6$\,:}\quad
Let $X_i=A_i$ for $i=1,2,3$ and $X_4=A_1\cup A_2$, $X_5=A_1\cup A_3$,
 $X_6=A_2\cup A_3$
The six graphs obtained in this way are $K_{a_1}$, $K_{a_2}$, $K_{a_3}$,
 $K_{a_1,a_2}$, $K_{a_1,a_3}$, $K_{a_2,a_3}$.
Each of them has average degree about $n/3$.

For an alternative, also here we cut both $B_1$ and $B_2$ into half:
 $B_1=B_1'\cup B_1''$ and $B_2=B_2'\cup B_2''$.
Then set $X_1=B_1$, $X_2=B_2$, $X_3=B_1'\cup B_2'$, $X_4=B_1''\cup B_2'$,
 $X_5=B_1'\cup B_2''$, $X_6=B_1''\cup B_2''$.
The average degrees are about $n/2$, $n/2$, $n/4$, $n/4$, $n/4$, $n/4$.

\msk

\nin
\underline{$k=7$\,:}\quad
This construction is obtained using the Fano plane, which is the first
 nontrivial particular case of Corollary \ref{c:blow-plane} $(i)$.
We will give all details there.
Currently  we do not know whether there is a substantially different
 alternative decomposition attaining $\Mad(G_1)+\ldots+\Mad(G_7)\sim 7n/3$.
\epf

\section{Exact formula for {\itshape M\/}$\ssl$({\itshape k}\hspace{1pt},\,{\itshape N}\hspace{1pt}) and sufficient condition for its equality with {\itshape M}\hspace{1pt}({\itshape k}\hspace{1pt},\,{\itshape n})}

\subsection{Tight relations: Steiner systems and graph designs}
\label{ss:main}

We begin with a complete determination of $M\ssl(k,N)$.

\btm
\label{t:list-k-n}
For any integers\/ $N\geq k\geq 2$, let\/ $p\geq 2$ be the unique integer
 for which\/ $k \binom{p}{2} \leq N < k \binom{p+1}{2}$ holds.
Write\/ $N$ in the form\/ $N = k \binom{p}{2} + qp + r$,
 where\/ $0\leq q < k$ and\/ $0\leq r < p$.
Then
$$
  M\ssl(k,N) = \begin{cases}
   \begin{tabular}{ll}
    $kp - k + q$ & if\/  $r\leq (p-1)/2$\,, \\
    $kp - k + q + 1 - 2(p-r)/(p+1)$ & if\/  $r\geq (p-1)/2$\,.
   \end{tabular}
    \end{cases}
$$
\etm

\bpf
Summing up the relevant values in the middle column of the following table
 (depending on $r$), the multiset $\{q*K_{p+1}, \, (k-q)*K_p\}$ for $r=0$
 or $\{q*K_{p+1}, \, (k-q-1)*K_p, \, 1*G_{p,r}\}$ for $0<r<p$
  shows that $M\ssl(k,N)$ is at least as large as given in the theorem.
We will prove that no list can have a larger $\Summad$ under the condition
 $\sum_{i=1}^k e(G_i) \leq N$.

\msk

\begin{center}
\begin{tabular}{lllll}
\quad graphs $H_i$ && total of $\Mad(H_i)$ && side condition \\ \hline
$q$ copies of $K_{p+1}$ && $qp$ && --- \\
$k-q-1$ copies of $K_p$ && $(k-q-1)(p-1)$ && --- \\
$1$ copy of $G_{p,r}$ && $p-1$ && $r\leq (p-1)/2$ \\
$1$ copy of $G_{p,r}$ && $p - 2(p-r)/(p+1)$ && $r > (p-1)/2$
\end{tabular}
\end{center}

\msk

Among the possibly many lists $\sL = \{G_1,\dots,G_k\}$ of $k$ graphs
 verifying $M\ssl(k,N)$, let us choose one satisfying the following
 hierarchic set of extremal conditions:
  \begin{itemize}
   \item[(0)]\quad $e(G_1)+\ldots+e(G_k)\leq N$ and $\summad = M\ssl(k,N)$.
   \item[(1)]\quad $\sL$ maximizes the number $s(\sL)$ of complete graphs,
    under (0).
   \item[(2)]\quad $\sL$ minimizes the difference $d(\sL)$ between the
    maximum order and the minimum order of complete graphs in the list,
    under (0) and (1).
   \item[(3)]\quad $\sL$ minimizes the number $h(\sL)$ of pairs that
    realize $d(\sL)$, under (0), (1), and (2).
  \end{itemize}
In addition we assume the following, which can be satisfied
 without loss of generality as we optimize a list for $M\ssl(k,N)$\,:
  \begin{itemize}
   \item[(4)]\quad Each $G_i\in\sL$ is of the form $G_i\cong G_{p_i,\,r_i}$
     with $0\leq r_i<p_i$, hence the representative of an $\cF_{p_i,\,r_i}$,
     having $\Mad(G_i)=g(m)$ where $m=e(G_i) = \binom{p_i}{2} + r_i$.
  \end{itemize}

We are going to prove that under these conditions there exists a list
 $\sL$ with $k-1 \leq s(\sL) \leq k$ and $d(\sL) \leq 1$.
Moreover, the possible unique non-complete member of $\sL$ is in tight
 relation with the $k-1$ complete members.

There are two types of non-complete graphs
 $G\in\fpr$ concerning $\Mad(G)$, namely
 type A: $0 < r \leq (p-1)/2$, and type B: $(p-1)/2 < r < p$.
However, if a $G \in \fpr$ of type A occurs in $\sL$, then
 it contains $r$ free edges.
Removing those edges, condition (0) is still respected, while the number
 $s(\sL)$ of complete members of $\sL$ is increased, a contradiction to (1).
Consequently, if there are any non-complete graphs in the list, all of
 them are of type B.

\msk

Suppose for a contradiction that $\sL$ contains two members of type B,
 say $G_i=G_{p_i,\,r_i}$ and $G_j=G_{p_j,\,r_j}$, and assume $p_i\leq p_j$.
We observe that $p_i=p_j$ must hold, for otherwise re-defining
 $G_i := G_{p_i,\,r_i+1}$ and $G_j := G_{p_j,\,r_j-1}$ would keep
 $e(G_i)+e(G_j)$ unchanged, while $\Mad(G_i)+\Mad(G_j)$ would increase
 with $\frac{2}{p_i+1} - \frac{2}{p_j+1} > 0$, contradicting (0).
So, let $p = p_i = p_j$, and suppose $p > r_j \geq r_i > (p-1)/2$.
We now re-define $G_j := K_{p+1}$ and $G_i := G_{p_i,\,r_i+r_j+1-p}$.
This keeps $e(G_i)+e(G_j)$ unchanged, does not decrease
 $\Mad(G_i)+\Mad(G_j)$ (in fact to avoid increase and hence a
  contradiction to (0) we must have $r_i+r_j\geq 3(p-1)/2$), but $G_j$
 becomes a complete graph, and so $s(\sL)$ increaes, contradicting (1).

\msk

Consequently, $k-1 \leq s(\sL) \leq k$ holds.
Next, we prove $d(\sL) \leq 1$.
Suppose for a contradiction that $G_i=K_{p_i}$ is a smallest clique in
 $\sL$, and $G_j=K_{p_j}$ is a largest clique in $\sL$, with $p_j>p_i+1$.
Let us re-define $G_i=K_{p_i+1}$ and $G_j=K_{p_j-1}$.
This keeps $\Mad(G_i)+\Mad(G_j)$ unchanged, and decreases $e(G_i)+e(G_j)$
 by $p_j-p_i-1>0$, so that (0) is respected.
Also, the number $s(\sL)$ of complete members in $\sL$ remains the same,
 respecting (1).
Should $G_i$ be unique smallest or $G_j$ unique largest in $\sL$, the
 modification would decrease $d(\sL)$, contradicting (2).
In any other case the value of $h(\sL)$ decreases, contradicting (3).

\msk

In this way, taking the fact $N<k\cdot\binom{p+1}{2}$ into account,
 one of the following lists is obtained:
 \begin{enumerate}
  \item \quad $\{k*K_{p'}\}$, with $p'\leq p$\,;
  \item \quad $\{q'*K_{p'+1}, \, (k-q')*K_{p'}\}$, with $p'\leq p$ and
   if $p'=p$ then $q'\leq q$\,;
  \item \quad $\{(k-1)*K_{p'}, \, 1*G_{p_k,\,r_k}\}$, with $p'\leq p+1	$\,;
  \item \quad $\{q'*K_{p'+1}, \, (k-q'-1)*K_{p'}, \, 1*G_{p_k,\,r_k}\}$.
 \end{enumerate}
The first and second cases obviously respect the formula stated in the
 theorem as an upper bound.
The proof for the third and fourth cases are completed by inspecting
 the relation between $p'$ and $p_k$ as it depends on whether $q'=0$
 or $q'>0$.

\msk

We claim that $q'=0$ implies $p_k=p'$ or $p_k=p'-1$.
Indeed, the assumption means that the first $k-1$ graphs
 all are copies of $K_{p'}$.
If $p_k\leq p'-2$, we move one edge from $G_1$ to $G_k$, hence re-define
 $G_1:=G_{p'-1,\,p'-2}$ and $G_k:=G_{p_k,\,r_k+1}$.
This keeps $e(G_1)+e(G_k)$ unchanged, but increases $\Mad(G_1)+\Mad(G_k)$
 with $\frac{2}{p_k+1} - \frac{2}{p'} > 0$, a contradiction to (0).
On the other hand, if $p_k>p'$,
 re-define $G_1:=K_{p'+1}$ and $G_k:=G_{p_k-1,\,r_k}$.
This decreases $e(G_1)+e(G_k)$ with $p_k-1-p'\geq 0$.
Further, $\Mad(G_1)$ is increased by 1, and $\Mad(G_k)$ is decreased
 from $\frac{p_k(p_k-1)+2r_k}{p_k+1} = p_k - 2 + \frac{2r_k+2}{p_k+1}$
  to $\frac{(p_k-1)(p_k-2)+2r_k}{p_k} = p_k - 3 + \frac{2r_k+2}{p_k}$.
Thus, the decrease concerning $G_k$ is
 $1 - 2(r_k+1)\bigl( \frac{1}{p_k} - \frac{1}{p_k+1} \bigr) < 1$,
 again a contradiction to (0).
As a consequence,
 we have $p'\leq p$ for $p_k=p'$, and $p'\leq p+1$ for $p_k=p'-1$,
 hence verifying the stated upper bound.

\msk

If $q'>0$, we show $p_k=p'\leq p$, and that $p'=p$ implies $q'\leq q$.
In this final case, copies of both $K_{p'}$ and $K_{p'+1}$ are present
 in the list $G_1,\dots,G_{k-1}$.
Referring to the proof of $d(\sL)\leq 1$ we see that $p_k<p'$ would contradict~(0)
 by moving an edge to $G_k$ from a copy of $K_{p'+1}$, whereas $p_k>p'$
 would yield a similar contradiction by removing a vertex of degree $p_k$
 from $G_k$ and extending a copy of $K_{p'}$ to $K_{p'+1}$.
Thus, here we must have $p'\leq p$, and the second claimed formula
 of the theorem applies.
\epf

We note that the following side product of this result also follows:

\begin{itemize}
 \item If $0<r\leq \frac{1}{2}\,(p-1)$, then $M\ssl(k,N)=M\ssl(k,N-j)$
   for $j=1,\dots,r$.
\end{itemize}
That is, in many cases it occurs that some (few or many) edges are
 unnecessary in the construction/verification of $M\ssl(k,N)$.

\btm
\label{t:Mad-k}
Let\/ $n\geq 3$ and\/ $2\leq k\leq \binom{n}{2}$ be given, and choose
  the integer\/ $p$ for which\/
  $k \binom{p}{2} \leq \binom{n}{2} < k \binom{p+1}{2}$ holds.
Let\/
 $q,r$ be the nonnegative uniquely determined integers such that\/
 $\binom{n}{2} = k \binom{p}{2} + qp + r$, where\/
 $0\leq q < k$ and\/ $0\leq r < p$. Then
$$
  M(k,n) \leq \begin{cases}
   \begin{tabular}{ll}
    $kp - k + q$ & if\/  $r\leq (p-1)/2$\,; \\
    $kp - k + q + 1 - 2(p-r)/(p+1)$ & if\/  $r\geq (p-1)/2$\,.
   \end{tabular}
    \end{cases}
$$
Moreover, equality holds whenever\/ $K_n$ admits an edge decomposition
 into
  \begin{itemize}
   \item $q$ copies of\/ $K_{p+1}$ and\/ $k-q$ copies of\/ $K_p$, if\/ $r=0$; or
   \item $q$ copies of\/ $K_{p+1}$, $k-q-1$ copies of\/ $K_p$,
    and one copy of any one graph\/ $G\in\fpr$, if\/ $0<r<p$.
  \end{itemize}
\etm

\bpf
Since $M(k,n) \leq M\ssl(k,\binom{n}{2})$, the claimed upper bound is a
 direct consequence of Theorem \ref{t:list-k-n}.
The maximum sum for lists has been calculated from a combination of
 complete graphs $K_p$ and $K_{p+1}$, possibly together with a
 representative graph $G_{p,r}\in\fpr$.
Thus, if $K_n$ admits an edge decomposition into the considered multiset
 of these graphs, then the list (multiset) also verifies that $M(k,n)$
 attains $M\ssl(k,\binom{n}{2})$.
\epf

Concerning the possible cases of equality, the
 condition given in the theorem is sufficient
 but not necessary in general.
This is shown by the following particular case.

\bex
Let\/ $k=7$ and\/ $n=8$.
We have\/ $e(K_8) = \binom{8}{2} = 28 = 12 + 12 + 4 = 2 e(K_4) + 4 e(K_3)
 + e(G_{3,1})$, hence\/ $M(7,8) \leq 2\, \Mad(K_4) + 4\, \Mad(K_3) +
  \Mad(G_{3,1}) = 6+8+2 = 16$.
Although two copies of\/ $K_4$ do not pack with five copies of\/ $K_3$ inside\/
 $K_8$, there exists a decomposition using the multiset\/ $\{ 4*K_3 ,
 2*(K_4-e) , 1*K_4\}$, which has\/ $\Summad = 4\cdot 2+2\cdot 2.5+3 = 16$.
Indeed, Let\/ $V(K_8)=\{v_1,\dots,v_8\}$.
Embed\/ $S(2,3,7)$ (the Fano plane) on the vertex set\/ $V(K_8)\smin\{v_8\}$
 in a way that the triplets\/ $(v_1,v_2,v_3)$, $(v_1,v_4,v_5)$,
 $(v_1,v_6,v_7)$ are blocks of the STS.
Extend these triplets with\/ $v_8$, to obtain\/ $K_4$ from the first and\/
 $K_4-e$ from the second and third (i.e., not using the edge\/ $v_1v_8$
 in the latter two).
The other four blocks of the STS remain unchanged, yielding the required
 copies of\/ $K_3$ for the decomposition of\/ $K_8$.
\eex

It is worth noting that the above construction was possible because
 the transformation $(K_4,G_{3,1}) \longleftrightarrow (K_4-e,K_4-e)$, where
 $G_{3,1}$ is actually the paw graph, keeps the sum of $\Mad$ values
  unchanged (while the individual values do change).

Further non-uniform extremal constructions will be presented later.
At this point we prove an asymptotic formula for $k$ growing
 proportionally to $\binom{n}{2}$.


\subsection{The upper range}
\label{ss:highest}

A Steiner system $S(2,p,n)$ is a collection of $p$-element
 subsets of an $n$-element underlying set, such that
 each pair of elements is contained in precisely one
 of those $p$-subsets.
An $S(2,p,n)$ can equivalently be interpreted as a set of
 $\binom{n}{2}/\binom{p}{2}$ mutually edge-disjoint
 copies of $K_p$ in $K_n$.
Beside $S(2,p,n)$ systems, there also exists
 a rich family of $\{K_p,K_{p+1}\}$-decompo\-sitions
 of $K_n$, and more generally $\{K_p,K_{p+1}\}$-packings.
In this way, many exact formulas can be derived
 from Theorem \ref{t:Mad-k}.
Below we explicitly formulate some of them.

\btm
\label{t:PSTS}
Let\/ $k=\binom{n}{2} - t$ with\/ $0\leq t\leq \frac{1}{3} (n-1)^2$.
 \begin{itemize}
  \item[$(i)$] If\/ $t$ is even, then\/ $M(k,n) = \binom{n}{2}-\frac{1}{2}t$,
   and the extremal multiset of subgraphs is\/
    $\{t/2 * K_3 , \, ( \, \binom{n}{2} - 3t/2 \, ) * K_2\}$.
  \item[$(ii)$] If\/ $t$ is odd, then\/
   $M(k,n) = \binom{n}{2}-\frac{1}{2}(t+1)+\frac{1}{3}$,
   and the extremal multiset of subgraphs is\/
    $\{(t-1)/2 * K_3 , \, ( \, \binom{n}{2} - 2 - 3(t-1)/2 \, ) * K_2 ,
     \, 1 * P_3\}$.
 \end{itemize}
\etm

\bpf
We apply Theorem \ref{t:Mad-k}.
Using its notation, the assumed range of $t$ yields $p=2$, hence every
 optimal edge decomposition can be transformed to the claimed multiset
 of copies of single edges and
 triangles, together with one copy of $P_3$ if $t$ is odd.
It is well known that such a packing of at least
 $\frac{1}{6} (n^2-2n-2)$ triangles exists indeed, for every $n$.

Suppose for a contradiction that there exists a further extremal
 multiset of graphs.
Consider the very last step of the process of transformations
 described in the proof of Theorem \ref{t:list-k-n}.
It is either $\{G_1,G_2\} \rightarrow \{2*K_3\}$ or
 $\{G_1,G_2\} \rightarrow \{K_3,P_3\}$.
Say, an edge $e$ is transposed from $G_1$ to $G_2$.
Then $G_1-e\cong K_3$, hence $G_1\cong K_3$+$\mathrm{leaf}$
 or $G_1\cong K_3\cup K_2$.
In either case, $e$ is a free edge in $G_1$.
However, $K_3$ and $P_3$ have no free edges, which means that the
 insertion of $e$ caused an increase of $\Mad(G_2)$.
This contradicts the assumption that we started with an extremal multiset.
\epf

\brm
The limitation for\/ $t$ in case of equality in Theorem \ref{t:PSTS}
 depends on the largest possible triangle packings in\/ $K_n$.
In the theory of block designs it was analyzed that the ``leave graph''
 (formed by the edges of\/ $K_n$ not covered by a largest partial Steiner
 triple system) is empty if\/ $n\equiv 1$ or\/ $3$ $(\mathrm{mod}$ $6)$,
 $\frac{n}{2} K_2$ if\/ $n\equiv 0$ or\/ $2$ $(\mathrm{mod}$ $6)$,
 $K_{1,3} \cup \frac{n-4}{2} K_2$ if\/ $n\equiv 4$ $(\mathrm{mod}$ $6)$,
 or just the\/ $4$-cycle $C_4$ if\/ $n\equiv 5$ $(\mathrm{mod}$ $6)$.
(See Theorem 3 of \cite{S-68}, or Table 40.22
 in Section VI.40.4 of \cite{CD-hand}.)
Accordingly, with a tight bound in Theorem \ref{t:PSTS}, $t$ can be
 as large as\/ $\frac{1}{3} n(n-1)$, $\frac{1}{3} n(n-2)$,
 $\frac{1}{3} (n-1)^2$, and\/ $\frac{1}{3} (n^2-n-5)$, respectively;
 because in the last two cases,
 one\/ $P_3$ can still be created optimally from\/ $K_{1,3}$ or\/ $C_4$.
\erm

\section{Lower-bound constructions via Steiner systems and graph designs}
\label{s:designs}

\subsection{Basic constructions}

\btm   \nev{Complete and partial Steiner systems}
\label{t:Stein}
Let\/ $p\geq 3$.
 \begin{itemize}
  \item[$(i)$] If\/ $\binom{p}{2} \mid \binom{n}{2}$ and\/ $p-1\mid n-1$,
   moreover if\/ $n>n_0(p)$,
   then for\/ $k=\frac{n(n-1)}{p(p-1)}$ we have\/ $M(k,n) = n(n-1)/p = (p-1)\cdot k$.

  \item[$(ii)$] As\/ $\nfty$, for $k = \floor{ \frac{n(n-1)}{p(p-1)} }$ we have
   $$
     \frac{n(n-1)}{p} - O(n) \leq M(k,n) \leq \frac{n(n-1)}{p} \,.
   $$
 \end{itemize}

\etm

\bpf
Under the given divisibility conditions of $(i)$, Wilson's theorem
 \cite{W-76} states that there exists an edge decomposition of $K_n$ into
 $\binom{n}{2}/\binom{p}{2}=k$ copies of $K_p$.
Each copy has $\Mad=p-1$, implying $M(k,n)\geq n(n-1)/p$.
This construction attains the maximum, due to
 Lemma \ref{t:conv} and Theorem~\ref{t:Mad-k}.

For $(ii)$ assume that $n$ is large enough with respect to $p$, and
 choose the largest integer $n'$ such that $n' \leq n$ and
 $n'\equiv 1$ or $p$ (mod $p(p-1))$; then $0 \leq n-n' \leq (p-1)^2$.
Observe that $n'-1$ is a multiple of $p-1$, and $\binom{n'}{2}$ is a
 multiple of $\binom{p}{2}$.
Hence, due to $(i)$, for $k'= \frac{n'(n'-1)}{p(p-1)}$ we have
 $M(k',n') = n'(n'-1)/p$.
Here $k-k' = \floor{\frac{n(n-1) - n'(n'-1)}{p(p-1)}}$, which obviously is
 smaller than $\binom{n}{2} - \binom{n'}{2} = \frac{n(n-1) - n'(n'-1)}{2}
 = (n-n')\,\frac{n+n'-1}{2} < \binom{p}{2} \cdot n$.
From the $\binom{n}{2} - \binom{n'}{2}$ edges of $K_n-K_{n'}$ we form
 $k-k'$ subgraphs of $K_n$ in an arbitrary way, and supplement them with
 a $K_p$-decomposition of $K_{n'}$.
Thus, we obtain $M(k,n) > M(k',n') = n'(n'-1)/p > n(n-1)/p - (p-1)\cdot n$.
\epf

\btm   \nev{Finite geometries}
\label{c:plane}
Let\/ $q$ be a prime power.
 \begin{itemize}
  \item[$(i)$] If\/ $n=k=q^2+q+1$, then\/ $M(k,n) = (q^2+q+1)q = qk =
   \frac{2}{q+1} \binom{n}{2}$.
  \item[$(ii)$] If\/ $n=q^2$ and\/ $k=q^2+q$, then\/ $M(k,n) = (q+1)q(q-1) =
   (q-1)k = \frac{2}{q} \binom{n}{2}$.
  \item[$(iii)$] If\/ $n=q^2+q$ and\/ $k=q^2+q+1$, then\/ $M(k,n) = q^3+q^2-1
    = qk - q - 1$.
%
  \item[$(iv)$] If\/ $n=q^2-1$ and\/ $k=q^2+q$, then\/ $M(k,n)
   = q^3 - 2q - 1 = (q-1)k - q - 1$.
  \item[$(v)$] If\/ $n=q^2-q$ and\/ $k=q^2+q-1$, then\/ $M(k,n)
   = q^3 - q^2 - 2 q + 1 = (q-2)\cdot k + q -1$.
 \end{itemize}
\etm

\bpf
Finite projective planes and in particular the projective Galois planes
 $PG(2,q)$, are Steiner systems $S(2,q+1,q^2+q+1)$, which are
 equivalent to decompositions of $K_n$ into copies of $K_{q+1}$.
Similarly, finite affine planes, including the affine Galois planes
 $AG(2,q)$, are Steiner systems $S(2,q,q^2)$, are equivalent to
  decompositions of $K_n$ into copies of $K_q$.
In this way the formulas if $(i)$ and $(ii)$ follow by
 Theorem \ref{t:Stein}.

For the other parts we make the following modifications:
 \begin{itemize}
  \item Concerning $(iii)$,
   pick a point of the projective plane $PG(2,q)$, and omit it from
   all the $q+1$ lines incident with it.
    This yields a set system over $q^2+q$ elements, with $q+1$ sets of
   size $q$ and $q^2$ sets of size $q+1$.
    The corresponding decomposition of $K_{q^2+q}$ has
     $$
       M(q^2+q+1,q^2+q) \geq \Summad = (q+1)(q-1) + q^3 = q^3+q^2-1 \,.
     $$
  \item Concerning $(iv)$,
   pick a point of the affine plane $AG(2,q)$, and omit it from
   all the $q+1$ lines incident with it.
    This yields a set system over $q^2-1$ elements, with $q+1$ sets of
   size $q-1$ and $q^2-1$ sets of size $q$.
    The corresponding decomposition of $K_{q^2-1}$ yields:

 \begin{eqnarray}
   M(q^2+q,q^2-1) & \geq & \Summad
    \nonumber \\
     & = & (q+1)(q-2) + (q^2-1)(q-1)
    \nonumber \\
     & = & (q+1)(q^2-q-1)
    \nonumber \\
     & = & q^3 - 2q - 1 \,.
    \nonumber 
 \end{eqnarray}
  \item Concerning $(v)$,
   pick a line $\ell$ of the affine plane $AG(2,q)$, and omit its
   points from all the $q^2$ lines which intersect $\ell$.
    This yields a set system over $q^2-q$ elements, with $q^2$ sets of
   size $q-1$ and $q-1$ sets of size $q$.
    (The $q-1$ lines parallel to $\ell$ remain intact.)
    The corresponding decomposition of $K_{q^2-q}$ has

 \begin{eqnarray}
   M(q^2+q,q^2-1) & \geq & \Summad \ = \
      q^2\cdot (q-2) + (q-1)^2
    \nonumber \\
     & = & q^3 - q^2 - 2 q + 1 \ = \
       (q-2)\cdot k + q -1 \,.
    \nonumber 
 \end{eqnarray}
 \end{itemize}
According to Theorem \ref{t:Stein}, all these three cases verify the exact
 value of $M(k,n)$, hence proving the formulas in $(iii)$--$(v)$.
\epf

The extremal structures in the cases of $(iii)$--$(v)$ above are
 non-uniform.
One can push this further, to obtain a larger class of tight
 constructions.
For this, we apply particular types of $K$-designs.
Given a set $K$ of positive integers, a pairwise balanced
 $K$-design ($K$-PBD, for short) of order $n$ and index $\lambda$
 is a pair $(X,\cB)$, where $X$ is an underlying set of cardinality $n$
 and $\cB$ is a collection of blocks (subsets of $X$) such that each
 pair $\{x,x'\}\subset X$ is contained in precisely $\lambda$ blocks.
A parallel class is a set of mutually disjoint blocks whose union is $X$.
A $K$-PBD $(X,\cB)$ is resolvable if $\cB$ can be partitioned into
 parallel classes.

\btm   \nev{Pairwise balanced K-designs with K = \{p, p{\bf \,+\,1}\}}

\ssk

\nin
$(i)$\quad
Let\/ $p\geq 3$ and\/ $K=\{p,p+1\}$.
Assume that\/ $(X,\cB)$ is a\/ $K$-PBD of order\/ $n$ and index\/ $1$,
 with all block sizes from\/ $K$
  (structurally equivalent to a\/ $\{K_p,K_{p+1}\}$-decomposition
 of\/ $K_n$), with\/ $b'$ blocks of size\/ $p$ and\/ $b''$ blocks of size\/ $p+1$.
Then, for\/ $k=b'+b''$, we have\/ $M(k,n) = (p-1)\cdot k + b''$.

\msk

\nin
$(ii)$\quad
Let\/ $p\geq 3$ and\/ $n\equiv p$ $(${\rm mod} $p(p-1))$.
Assume that at least one of the following conditions holds:
 \begin{itemize}
  \item  $n > \exp\{\exp\{p^{12p^2}\}\}$
   and\/ $n'\leq \frac{n-1}{p-1}$;
  \item  $n' < 10^{-4} n / p^{5/2}$.
 \end{itemize}
Then the following holds.
 \begin{itemize}
  \item[$(a)$] If there exists an\/ $S(2,p,n')$-system, then for\/
   $k= \frac{n(n-1)}{p(p-1)} + \frac{n'(n'-1)}{p(p-1)}$ we have\/
    $M(k,n+n') = (p-1)\cdot k + (n\cdot n')/p$.
  \item[$(b)$] If there exists an\/ $S(2,p+1,n')$-system, then for\/
   $k= \frac{n(n-1)}{p(p-1)} + \frac{n'(n'-1)}{p(p+1)}$ we have\/
    $M(k,n+n') = (p-1)\cdot k + (n\cdot n')/p + n'(n'-1)/(p+1)$.
 \end{itemize}
\etm

\bpf
Similarly as above, Theorem \ref{t:Stein} implies that every decomposition
 of $K_n$ into any $b'$ copies of $K_p$ and $b''$ copies of $K_{p+1}$
 provides the exact value of $M(b'+b'',n)$.
This verifies $(i)$.

In order to prove $(ii)$ we need to show that the corresponding
 decompositions of $K_{n+n'}$ exist under the conditions imposed on
 $n$, $n'$, and $p$.
We view $K_{n+n'}$ as the complete join $K_n+K_{n'}$.

For the case of $n' \leq \frac{n-1}{p-1}$ and
 $n > \exp\{\exp\{p^{12p^2}\}\}$, Chang \cite[Theorem~7.2]{C-00} proved
 that there exists a resolvable $S(2, p, n)$ system.
Using its blocks as vertex sets of complete subgraphs in $K_n$,
 from the parallel classes of a resolution we obtain a decomposition
 of $K_n$ into $\frac{n-1}{p-1}$ copies of $\frac{n}{p} K_p$.

Each of the $n'$ vertices in $K_{n'}$ can be used to extend the $n/p$
 copies of $K_p$ in a parallel class to the same number of $K_{p+1}$
 subgraphs in $K_{n+n'}$, distinct classes taken for distinct vertices.
Note that these cover all the edges connecting $V(K_n)$ to $V(K_{n'})$.
The other $\frac{n-1}{p-1} - n'$ parallel classes and the blocks of the assumed
 $S(2,p,n')$ or $S(2,p+1,n')$ system ensure that a required decomposition
 of $K_{n+n'}$ into $k$ complete subgraphs (each of them being
 $K_p$ or $K_{p+1}$) exists, indeed.

For the case of $n'$ much smaller than $n/p$ we need not assume a huge
 $n$; namely, the condition $p\mid n$ will suffice.
Also here we shall need $n'$ parallel classes $\frac{n}{p} K_p$, which
 we find using the Hajnal--Szemer\'edi Theorem \cite{HSz-70}.
It states that, for every positive integer $s$, every graph with maximum
 degree at most $s$ has an equitable coloring with $s + 1$ colors.
Putting $s=n/p - 1$, equitable means that all color classes have size $p$;
 those classes will be vertex sets of the copies of $K_p$, together
 forming parallel classes $F_i\cong \frac{n}{p} K_p$.
We select $F_1,F_2,\dots,F_{n'}$ sequentially.
The first one, $F_1$, is just any choice of $\frac{n}{p} K_p$.
Let us denote $E_i=E(F_i)$.
Having chosen $F_1,\dots,F_j$, consider the graph whose edge set is
 $E_1 \cup \cdots \cup E_j$.
This graph is regular of degree $(p-1)j<n/p-1$, therefore it admits an
 equipartite coloring with $n/p$ colors.
This coloring yields $F_{j+1}$.
The procedure leads to $n'$ parallel classes as needed.
The edge set of each $F_i$ is supplemented with a distinct vertex of
 $K_{n'}$, hence yielding $(n\cdot n')/p$ copies of $K_{p+1}$.
We can also use the blocks of the assumed
 $S(2,p,n')$ or $S(2,p+1,n')$ to obtain further complete subgraphs
 $K_p$ or $K_{p+1}$ in $K_{n+n'}$.

It remains to decompose the graph $H$ having the edge set
 $E(H) = E(K_n)-\cup_{i=1}^{n'} E_i$ into copies of $K_p$.
Gustavsson \cite{G-91} proved that this can be done whenever the following
 conditions are satisfied: $\binom{p}{2} \mid e(H)$, the degree of every
 vertex is a multiple of $p-1$, and $\delta(H)>(1-\epsilon)n$, where
 $\epsilon=\epsilon(p)>0$ is a sufficiently small absolute constant.
Making a substantial improvement on $\epsilon$, Glock et al.\
 \cite[Corollary 1.6]{GKLMO} proved that $\epsilon \geq 10^{-4}p^{-3/2}$
 is a suitable choice for $K_p$-decomposition.
Thus, the described strategy is feasible whenever at most
 $\epsilon n / (p-1)$ parallel classes are needed.
This completes the proof of the theorem.
\epf

\brm
The principle of proof for part\/ $(ii)$ also yields the following more
 general recursion:
\begin{itemize}
 \item[$(iii)$] Let\/ $n,n'\in\nnn$, $K=\{p,p+1\}$, and suppose that there
  exist block designs of the following types: a PBD$(K,n)$ system\/
   $(X_n,\cB_n)$ with\/
  $b'_n$ blocks of size\/ $p$ admitting at least\/ $n'$ parallel classes
   and with\/ $b''_n$ blocks of size\/ $p+1$ ($b_n'+b_n''=|\cB_n|$),
    and a PBD$(K,n')$ system\/ $(X_{n'},\cB_{n'})$ with\/
  $b'_{n'}$ blocks of size\/ $p$ and with\/ $b''_{n'}$ blocks of size\/ $p+1$
   ($b_{n'}'+b_{n'}''=|\cB_{n'}|$).
   Then, for\/ $k = |\cB_n| + |\cB_{n'}| + (n\cdot n')/p$ we have
  $$
    M(k, n+n') = (p-1)\cdot k + b_n'' + b_{n'}'' + (n\cdot n')/p \,.
  $$
\end{itemize}
\erm

It was shown in Theorem \ref{t:PSTS} that the behavior of $M(k,n)$
 for $k\geq \frac{1}{3} \binom{n}{2}$ is understood fairly well.
The next result gives some information concerning the remaining part of
 the quadratic range of $k$.

\btm   \nev{Fixed proportional growth of {\bf\itshape k}}
\label{t:proport}
Let\/ $r\geq 3$ be any rational number.
There exist infinitely many\/ $n$ such that\/ $M(k,n)$
 for\/ $k = \frac{1}{r} \binom{n}{2}$ is determined by a\/
 $\{K_p,K_{p+1}\}$-design, where\/ $p\geq 3$ is the integer satisfying\/
 $\binom{p}{2} \leq r < \binom{p+1}{2}$.
\etm

\bpf
Let us note first that the case $r=\binom{p}{2}$ is already settled in
 part $(i)$ of Theorem \ref{t:Stein}.
Therefore, in the sequel we assume $\binom{p}{2} < r < \binom{p+1}{2}$.

We search for an edge decomposition of $K_n$ into the members of the
 multiset $\textstyle\{xk * K_p \, , (1-x)k * K_{p+1}\}$,
 with a suitable choice of $x\in(0,1)$.
Since $e(K_n) = \binom{n}{2} = rk$, this requires
 $$
   xk\binom{p}{2} + (1-x)k \binom{p+1}{2} = rk \,,
 $$
 thus
 $$
   x = \frac{p+1}{2} - \frac{r}{p} \,.
 $$
Since the integer $p$ is determnimned by $r$, we see that $x$ is a
 function of $r$, a rational number between 0 and 1 because
 $\frac{p-1}{2} < \frac{r}{p} < \frac{p+1}{2}$ by the choice of~$p$.
In particular, $x$ is independent of the actual choice of the number $n$
 of vertices.

Assume $x=s/t$ ($s,t\in\nnn$\,, \,$1\leq s\leq t-1$),
 and let $H = sK_p \cup (t-s)K_{p+1}$.
Then $\delta(H)=p$ and $\Delta(H)=p+1$.
Thus, Wilson's theorem implies that $K_n$ admits an $H$-decomposition
 whenever $\binom{n}{2}$ is a multiple of $e(H)$ and $n$ is
 sufficiently large.
Here $e(H)=rt$, hence $K_n$ is decomposed into $\binom{n}{2}/(rt)$
 copies of $H$, which then further decompose into
 $\frac{s}{rt}\binom{n}{2}$ copies of $K_p$ and
 $\frac{t-s}{rt}\binom{n}{2}$ copies of $K_{p+1}$.
Consequently, a $\{K_p,K_{p+1}\}$-design is obtained with exactly
 $\frac{1}{r}\binom{n}{2}$ subgraphs.
By the ``moreover'' part of Theorem \ref{t:Mad-k} we obtain:
 $$
   M(k,n) = \left( \frac{s\cdot(p-1)}{rt} +
    \frac{(t-s)\cdot p}{rt} \right) \! \binom{n}{2}
     = \frac{tp-s}{rt} \binom{n}{2} ,
 $$
  implying the assertion.
\epf

In particular, we obtain:

\bcr
For any rational\/ $r\geq 1$ there are infinitely many\/ $n$ such that\/
 $M(k,n)$ is an integer for\/ $k = \frac{1}{r} \binom{n}{2}$.
\ecr

\bpf
Assuming $r=a/b$, apply the formula at the end of the proof of
 Theorem \ref{t:proport}, and restrict the solutions to
 those $n$ for which $\binom{n}{2}$ is a multiple of $ta$.
\epf

\subsection{Constructions via substitution and the asymptotic growth of {\itshape M}\hspace{1pt}({\itshape k}\hspace{1pt},\,{\itshape n})}
\label{ss:triang+stein}

Here we first consider the cases of triangular numbers $k=\binom{t+1}{2}$.
Since $M(k,n) < \sqrt{k} \, n$ by Corollary~\ref{c:M-sqrt-k}, the
 following lower bound is tight apart from a multiplicative constant
 in the entire range of $k$.
Substantial improvement for large $k$ will be given later in this section.

\bpn
Let\/ $k=\binom{t+1}{2}$ for some\/ $t\geq 2$, and assume\/ $t\mid n$.
Then\/ $M(k,n) \geq \frac{t+1}{2}\,n - t > \sqrt{k/2} \, (n - 2)$.
\epn

\bpf
Partition $V(K_n)$ into $t$ sets $V_1,\dots,V_t$ of cardinality $n/t$ each.
Define $t+\binom{t}{2}$ subgraphs indexed as $G_i$ for $1\leq i\leq t$
 and $G_{i,j}$ for $1\leq i<j\leq t$, as follows.
Each $G_i$ is a copy of $K_{n/t}$ with $V(G_i)=V_i$, and each $G_{i,j}$
 is a copy of $K_{\frac{n}{t},\frac{n}{t}}$ with $V(G_{i,j})=V_i\cup V_j$.
All these graphs are regular.
More explicitly, we have obtained $t$ graphs of degree $n/t - 1$, and
 $\frac{1}{2} t(t-1)$ graphs of degree~$n/t$.
Consequently, $M(k,n) \geq n - t + \frac{1}{2} (t-1) n $.
The last inequality is valid because $k<(t+1)^2\!/2$.
\epf

\btm
\label{t:blow-Stein}
Let\/ $v>r>2$ be such that a Steiner system\/ $S(2,r,v)$ exists.
If\/ $k=\frac{v(v-1)}{r(r-1)}$, and\/ $n$ is any integer divisible by\/ $v$,
 then\/ $M(k,n)\geq \frac{v}{r}\left(n-1\right)$.
\etm

\bpf
Let $(X,\cB)$ be an $S(2,r,v)$ with point set $X$ and blocks
 $B_1,\dots,B_k$.
Replace the points $x\in X$ with mutually disjoint sets $V_x$ of
 cardinality $n/v$.
Then, viewing $\cup_{x\in X} V_x$ as the vertex set of $K_n$, we obtain
 $k$ sets $Y_i := \cup_{x\in B_i} V_x$ of cardinality $rn/v$.
For $i=1,\dots,k$ let $H_i$ be the complete $r$-partite graph whose edges
 are the vertex pairs $y,y'$ where $y,y'\in Y_i$, say $y\in V_x$ and
 $y'\in V_{x'}$, such that $x\neq x'$.
Distribute the remaining edges of $K_n$ (i.e., those which are contained
 in some $V_x$) arbitrarily among the $H_i$.
In this way we obtain an edge decomposition of $K_n$ into $k$ subgraphs;
 denote them by $G_1,\dots,G_k$.

The average number of edges in the $G_i$ is $\frac{n(n-1)}{2k}$, and
 each of them has order $rn/v$; hence the average of their $d(G_i)$ is 
 equal to $\frac{v}{rk}\left(n-1\right)$.
Thus, $M(k,n) \geq \Summad \geq \sum_{i=1}^k d(G_i) =
 \frac{v}{r}\left(n-1\right)$.
\epf



\bcr
\label{c:blow-plane}
Let\/ $q$ be a prime power.
 \begin{itemize}
  \item[$(i)$] If\/ $k=q^2+q+1$ and\/ $k\mid n$, then\/ $M(k,n) \geq
   (q + \frac{1}{q+1}) (n-1)$.
  \item[$(ii)$] If\/ $k=q^2+q$ and\/ $q^2\mid n$, then\/ $M(k,n) \geq qn - q$.
 \end{itemize}
\ecr

\bpf
Apply Theorem \ref{t:blow-Stein} to the projective planes
 $S(2,q+1,q^2+q+1)$ and affine planes $S(2,q,q^2)$, with $v=q^2+q+1$
 or $v=q^2$ and $r=q+1$ or $r=q$, respectively.
\epf

Now we are in a position to prove that $M(k,n)$ as a function of $n$
 grows essentially with $\sqrt{k}\,n$.

\btm
\label{t:sqrt-k}
For every\/ $n$ and every fixed\/ $k$ we have
 $$
   ( 1 - \epsilon_k ) \sqrt{k}\, n < M(k,n) < \sqrt{k}\, n \,,
 $$
  where\/ $\epsilon_k\to 0$ as\/ $k\to\infty$.
In convergence,
  $$
    \lim_{k\to\infty}
   \Bigl( \liminf_{n\to\infty} \frac{M(k,n)}{\sqrt{k}\,n} \Bigr) =
    \lim_{k\to\infty}
   \Bigl( \limsup_{n\to\infty} \frac{M(k,n)}{\sqrt{k}\,n} \Bigr) = 1 \,.
  $$
\etm

\bpf
Recall that the upper bound $\sqrt{k}\,n$ is valid by Corollary \ref{c:M-sqrt-k}.
Moreover, if $k = q^2+q+1$ where $q$ is a prime power and $n$ is a
 multiple of $k$, then the lover bound holds in a stronger form.
Namely, we can write
   $$(\sqrt{k} - 1/2)\,n - c_k$$
 instead of
 $( 1 - \epsilon_k ) \sqrt{k}\,n$, choosing $q + \frac{1}{q+1}$
 (which is less than $\sqrt{k}$\,) for the constant $c_k$.
Indeed, from part~$(i)$ of Corollary~\ref{c:blow-plane} we obtain
 $$
   M(k,n) + q + \frac{1}{q+1} \geq ( q + \frac{1}{q+1} ) \, n >
   ( \sqrt{q^2+q+1} - 1/2) \, n = (\sqrt{k} - 1/2)\,n
 $$
  where the inequality in the middle can easily be checked.

For orders not divisible by $k$, let us consider $n'=n+t=ks+t$, where
 $1\leq t\leq k-1$.
By monotonicity with respect to $n$ (due to Proposition \ref{p:monot}) we know
\begin{eqnarray}
  M(k,n') & \geq & M(k,n) \ > \ (\sqrt{k} - 1/2)\,(n'-t) - \sqrt{k}
   \nonumber \\
  & \geq & (\sqrt{k} - 1/2)\,n' - \sqrt{k} - (k-1)(\sqrt{k} - 1/2)
   \nonumber \\
  & = & (\sqrt{k} - 1/2)\,n' - k^{3/2} + k/2 - 1/2 \,.
   \nonumber
\end{eqnarray}
Hence $k^{3/2}$ is a suitable choice for $c_k$
 whenever $k=q^2+q+1$ with any prime power $q$.

Suppose that $k$ cannot be written in this form.
Let then $q$ be the smallest prime power such that $q^2+q+1>k$.
It is well known that there exists a sublinear function $f(x)=o(x)$ such
 that the interval $[x-f(x),x]$ contains a prime number for every $x>2$.
(The current record concerning the growth of $f$ seems to be $o(x^{0.525})$,
 see \cite{BHP-01}.)
So, we can choose a prime (or prime power) $q'$ with $q-q'=o(q)$.
Let $k' = (q')^2 + q' +1$; certainly $k'<k$, moreover
 $k-k' = o(k)$ and $\sqrt{k}-\sqrt{k'} = o(\sqrt{k}\,)$ as $k$ grows.
Since $M(k,n)$ is monotone increasing in $k$, again by
 Proposition \ref{p:monot}, we obtain:
 $$
   M(k,n) > M(k',n) > (\sqrt{k'} - 1/2)\,n - c_{k'} \geq
     (\sqrt{k} - \epsilon_k \sqrt{k}\,)\,n
 $$
  with a suitably chosen $\epsilon_k>0$.

For the limits $s_k := \liminf_{n\to\infty} \frac{M(k,n)}{\sqrt{k}\,n}$
 and $t_k := \limsup_{n\to\infty} \frac{M(k,n)}{\sqrt{k}\,n}$ this yields
 $1-o_k(1) \leq s_k \leq t_k \leq 1$, implying the theorem.
\epf

In the last inequality above, the constant $\epsilon _k$ can be estimated
 via the $f(x) = o(x^{0.525})$ theorem quoted in the proof.
The extension of the theorem for non-fixed $k$ faces the difficulty that
 the $o$-bound on $f$ requires $x$ to be large, but $n$ puts limitations
 on that.

Part $(i)$ of Corollary \ref{c:blow-plane} can be extended also for
 values of $k$ which are not of the form $q^2+q+1$.
This offers a replacement of the above error term
 $\epsilon _k \sqrt{k}\,n$ with an explicit function.

\bpn
\label{p:plane+r}
Let\/ $k=q^2+q+1+r$, where\/ $q$ is a prime power and\/ $1\leq r\leq q^2+q+1$.
If\/ $(q^2+q+1) | n$, then\/ $M(k,n) \geq (q + \frac{1}{q+1})\,n +
 \frac{rq}{(q+1)(q^2+q+1)}\,n - (q+r)$.
\epn

\bpf
For this purpose we need to specify the construction of
 Theorem \ref{t:blow-Stein} more strictly.
We apply the fact that every projective Galois plane has a
 cyclic representation (e.g., via difference sets).
Let $x_1,\dots,x_{q^2+q+1}$ and $\ell_1,\dots,\ell_{q^2+q+1}$ be the
 points and lines, respectively, of $PG(2,q)$, such that $x_i\in\ell_i$
 for $i=1,\dots,q^2+q+1$ and the rotation $x_i\mapsto x_{i+1}$ performed
 modulo $q^2+q+1$ is an automorphism of $PG(2,q)$.
Replace each $x_i$ with a set $X_i$ of cardinality $n/(q^2+q+1)$, creating
 a complete equipartite graph $H_i$ from $\ell_i$ with $q+1$ vertex classes.

If $r+1\leq i \leq q^2+q+1$, we let $G_i$ be the graph obtained from $H_i$
 by inserting the edges of the complete graph whose vertex set is $X_i$.
This $G_i$ has $\frac{1}{q+1}|G_i|$ vertices of degree $|G_i|-1$ (namely the
 vertices in $X_i$), and its other $(1-\frac{1}{q+1})|G_i|$ vertices have
 degree $(1-\frac{1}{q+1})|G_i|$.
Recalling that $|G_i| = \frac{q+1}{q^2+q+1}\,n$, we obtain
 \begin{eqnarray}
  \sum_{i=r+1}^{q^2+q+1}
   d(G_i) & = &
  \sum_{i=r+1}^{q^2+q+1} \frac{1}{|G_i|} \left( \frac{1}{q+1}|G_i|
   \left( |G_i|-1 \right) + \left( 1-\frac{1}{q+1} \right)^2|G_i|^2 \right)
     \nonumber \\
   & = &
  \sum_{i=r+1}^{q^2+q+1} \left( \left( \frac{1}{q+1} + \left( \frac{q}{q+1} \right)^2 \right) |G_i|
    - \frac{1}{q+1} \right)
     \nonumber \\
   & = & \left( q^2+q+1-r \right) \left( \frac{q^2+q+1}{(q+1)^2} \left( \frac{q+1}{q^2+q+1} \, n \right) - \frac{1}{q+1} \right)
     \nonumber \\
   & \geq &
     \frac{q^2+q+1-r}{q+1}\,n - q \,.   \label{egy}
 \end{eqnarray}

If $1\leq i\leq r$, we simply put $G_i=H_i$.
These graphs are regular of degree $(1-\frac{1}{q+1})|G_i|$, hence
 \begin{eqnarray}
  \sum_{i=1}^{r}
   d(G_i) & = & r \left( \left( 1 -\frac{1}{q+1} \right) |G_i| \right)
     \nonumber \\
   & = & \frac{rq}{q+1} \, \frac{q+1}{q^2+q+1} \,n \,.   \label{ket}
 \end{eqnarray}

Finally, if $q^2+q+2\leq i\leq q^2+q+1+r$, we let $G_i$ be the complete
 graph on the vertex set $X_i$.
These $r$ graphs are regular of degree $|X_i|-1$, hence
 \begin{eqnarray}
  \sum_{i=q^2+q+2}^{q^2+q+1+r}
   d(G_i) & = & \frac{r}{q^2+q+1} \,n - r .   \label{har}
 \end{eqnarray}

Summing up (\ref{egy})+(\ref{ket})+(\ref{har}) we obtain
 \begin{eqnarray}
  M(k,n) & \geq & \Summad 
     \nonumber \\
     & = & \left( q + \frac{1}{q+1} \right)\,n +
      \left( \frac{q+1}{q^2+q+1} - \frac{1}{q+1} \right)\, rn - q - r \,,
     \nonumber
 \end{eqnarray}
  proving the assertion.
\epf

An explicit lower bound of similar flavor can also be given for the cases
 where $k$ is big, but
 not much bigger than the number of blocks in a Steiner system.

\bpn
Assume that\/ $k=\binom{n}{2}/\binom{p}{2}$, and that a Steiner system\/
 $S(2,p,n)$ exists.
If\/ $1\leq r\leq 2k$, then\/ $M(k+r,n) \geq
 M(k,n) + rp/4 - cr$ for a constant\/ $c$.
\epn

\bpf
We know from Theorem \ref{t:Mad-k} that $M(k,n)$ for $k=\binom{n}{2}/\binom{p}{2}$
 is generated via $S(2,p,n)$ whenever such a Steiner system exists.
In that construction, every graph $G_i$ in the decomposition of $K_n$ is a
 copy of $K_p$.
Applying Theorem~\ref{t:k=2} and the first subcase of
 Proposition \ref{p:3--7}, we decompose selected copies of $K_p$ into
 two subgraphs $H_{i,1},H_{i,2}$ with
 $\Mad(H_{i,1})+\Mad(H_{i,2})\geq 5p/4-3/2$
 or into three subgraphs $H_{i,1},H_{i,2},H_{i,3}$ with
 $\Mad(H_{i,1})+\Mad(H_{i,2})+\Mad(H_{i,3})\geq 3p/2-c_3$.
Increasing the number of subgraphs with $r$, this operation increases
 $\summad$ with at least $rp/4-cr$, where $c$ is relatively small.
\epf

Also, the following recursive lower bound can be stated.

\bpn
For every\/ $a,b,t\in\nnn$ with\/ $b\geq 3$, $t\geq 2$ and\/
 $1\leq a\leq \binom{b}{2}$ we have\/ $M(a+b,tb) \geq t\cdot(M(a,b)+b)-b$.
\epn

\bpf
Applying the principle of the preceding constructions, replace each vertex
 $x_i$ of a structure verifíying $M(a,b)$ with a set $X_i$ of cardinality
 $t$, blow up the corresponding graphs $G_1,\dots,G_a$ accordingly
 (hence multiplying their $\Mad$ values by $t$), and create $b$ further
 graphs $G_{a+1},\dots,G_{a+b}$ regular of degree $t-1$, $G_{a+i}$ being
 the complete graph on the vertex set $X_i$ for $i=1,\dots,b$.
\epf

Applying the above implication for $t=n/b$, we may abbreviate the
 conclusion as $M(a,b) \Longrightarrow M(a+b,n) \gtrsim c_{a,b} \cdot n$
 for an asymptotic lower bound as $n$ gets large,
 where $c_{a,b}$ abbreviates the value $\frac{1}{b} \cdot M(a,b) + 1$.
Some pairs $a,b$ provide lower bounds that coincide with those given in
 Proposition \ref{p:3--7}:
  \begin{itemize}
   \item $M(1,2) \Longrightarrow M(3,n) \gtrsim 3n/2$;
   \item $M(1,3) \Longrightarrow M(4,n) \gtrsim 5n/3$;
   \item $M(3,3) \Longrightarrow M(6,n) \gtrsim 2n$.
  \end{itemize}
In the same way also $M(k,n) \Longrightarrow
 M(k+n,tn) \geq t\cdot M(k,n) + (t-1)\cdot n$ can be obtained.

\section{Applications to other graph invariants}
\label{s:appli}

In this section we derive implications of the above results for other
 graph invariants, namely:
 \begin{itemize}
  \item clique number, $\omega(G)$ (largest number of mutually adjacent
   vertices in $G$);
  \item chromatic number, $\chi(G)$ (minimum number of colors in a vertex
   coloring such that no two adjacent vertices have the same color);
  \item choice number or list chromatic number, $\ch(G)$ (minimum size $m$ of
    color lists $L_v$ assigned to the vertices $v\in V(G)$ guaranteeing that
    there exists a coloring $\vp:V(G)\to\cup_{v\in V(G)} L_v$ with $\vp(v)\in L_v$
    for all $v\in V(G)$ and $\vp(v)\neq \vp(w)$ for all $vw\in E(G)$
    whenever $|L_v|=m$ for all $v\in V(G)$);
  \item degeneracy, $\dg(G)$ (smallest integer $d$ such that $G$ admits a
   vertex order $v_1,\dots,v_n$ in which each $v_j$ has at most $d$
   preceding neighbors $v_i$ with $i<j$, equal to $\max\{\delta(H)\}$
   taken over all induced subgraphs $H\subseteq G$);
  \item coloring number, $\col(G)=\dg(G)+1$, also called Szekeres--Wilf number;
  \item highest induced vertex connection, $\kaps(G)$ (maximum vertex
   connection number $\kappa(H)$ over all induced subgraphs $H\subseteq G$);
  \item highest induced edge connection, $\lams(G)$ (maximum edge
   connection number $\lambda(H)$ over all induced subgraphs $H\subseteq G$).
 \end{itemize}
Applying the notational convention described in the introduction, we write
 $\omega(k,n)$, $\chi(k,n)$, $\dg(k,n)$, $\col(k,n)$, $\kaps(k,n)$, $\lams(k,n)$
 for the largest possible sum $\sum_{j=1}^k \pp(G_j)$ taken over all
 edge decompositions of $K_n$ into $k$ subgraphs $G_1,\dots,G_k$, where
 $\pp\in\{\omega, \chi, \dg, \col, \kaps, \lams\}$.

Since the inequality chains
 $\Mad(G)+1\geq\col(G)\geq\ch(G)\geq\chi(G)\geq\omega(G)$ and
 $\Mad(G)\geq\dg(G)\geq\lams(G)\geq\kaps(G)$
 are valid for every graph $G$, also the corresponding extremal functions
 satisfy the following numerical relations
 for all $n\geq 3$ and $2\leq k\leq \binom{n}{2}$\,:
  $$
    \omega(k,n) \leq \chi(k,n) \leq \ch(k,n) \leq \col(k,n) \leq M(k,n) + k \,,
  $$
  $$
    \kaps(k,n) \leq \lams(k,n) \leq \dg(k,n) \leq M(k,n) \,.
  $$
Bickle \cite{Bic-12} conjectures that
 $\dg(k,n) = \kaps(k,n) = \lams(k,n)$ holds for all $k,n$.
We shall see in Theorem \ref{t:pp-k-n} that this is indeed the case
 at least when a supplementary assumption is imposed.

F\"uredi et al.\ studied four of the above seven parameters; next we
 summarize their findings.

\btm   [F\"uredi, Kostochka, Stiebitz, \v Skrekovski, West \cite{FKSSW}]
\label{t:F-etal}
\
  \begin{enumerate}
   \item 
    \ $\omega(k,n) \leq n+\binom{k}{2}$,
    with equality for all\/ $n\geq \binom{k}{2}$.
   \item 
    \ $\chi(k,n) \leq n+7^k$.
   \item 
    \ $\ch(k,n) \leq n+3k!\sqrt{1+8n\ln n}$,
    moreover\/ $\ch(k,n) \geq n+ck\ln(n/\!\sqrt{k}\,)$ for a constant\/ $c>0$
    if\/ $k$ is a triangular number, $k=\binom{t+1}{2}$, and\/ $t\mid n$.
   \item 
    \ $\col(k,n) \leq \sqrt{k}\,n+k$,
    moreover\/ $\col(k,n) \geq (\sqrt{k}-1)n+k$
    if\/ $k=q^2+q+1$ for a prime power\/ $q$, and $k\mid n$.
   \item 
    \ $\col(k,n) \leq (k-1)(n+1)/2$
    if\/ $k\geq 5$; for smaller\/ $k$,
    $\col(3,n) = (3n+3)/2$ and\/ $\col(4,n) = (5n+7)/3$.
  \end{enumerate}
\etm

Let us recall further the following Norhaus--Gaddum-type result
 (converted to the present notation) from the
 paper which initiated the study of choice number.

\btm [Erd\H os, Rubin, Taylor \cite{ERT}]
\label{t:ERT}
For every\/ $n\in\nnn$,
$\ch(2,n) \leq \col(2,n) \leq n+1$.
\etm

Comparing Theorem \ref{t:ERT} with Theorem \ref{t:k=2} it turns out that
 $M(2,n)\sim \frac{5}{4}\,n$ is much bigger than $\col(2,n)$.
For $k=3$ and $k=4$ it is not clear whether the situation
 between $M(k,n)$ and $\col(k,m)$ is similar,
 as the corresponding constructions in Proposition \ref{p:3--7} yield
 lower bounds which coincide with those in Part 5 of Theorem \ref{t:F-etal}.
However, as $k$ gets large and $n\to\infty$, the relative gap between
 $\col(k,n)$ and $M(k,n)$ gets small.
What is more, our results on $M(k,n)$ improve substantially the earlier
 estimates on several functions $\pp(k,n)$ in certain ranges of the
 parameters.

Let us recall that the list version $M\ssl(k,\binom{n}{2})$ of $M(k,n)$
 (and more generally, $M\ssl(k,N)$ for all $N\geq k$) was
 determined exactly for all $n\geq 3$ and all $2\leq k\leq \binom{n}{2}$
 in Theorem \ref{t:list-k-n}.
It has also been proved that every list of graphs being extremal for
 $M\ssl(k,\binom{n}{2})$ can be transformed to one of two types:
  \begin{itemize}
   \item Type 1: $q$ copies of $K_{p+1}$ and $k-q$ copies of $K_p$, if $r=0$;
   \item Type 2: $q$ copies of $K_{p+1}$, $k-q-1$ copies of $K_p$,
    and one copy of any one graph $G\in\fpr$, if $0<r<p$;
  \end{itemize}
 where $\binom{n}{2} = k \binom{p}{2} + qp + r$, with
 $0\leq q < k$ and $0\leq r < p$.
As stated in Theorem \ref{t:Mad-k}, we have $M(k,n) = M\ssl(k,\binom{n}{2})$
 whenever $M(k,n)$ is attained by a packing of Type 1 or Type 2.
Moreover, a bunch of results above show that this situation occurs
 for infinitely many pairs $k,n$.
Consequences can be stated for other parameters as well.

\btm
\label{t:pp-k-n}
If\/ $M(k,n)$ is attained by a packing of Type 1 or Type 2, and in case of
 Type 2 with\/ $r\geq (p-1)/2$ the packing is established using a graph\/
  $G\in \cG_{p,r}$ with\/ $\omega(G)=p$, then
 $$
   \omega(k,n) = \chi(k,n) = \ch(k,n) = \col(k,n) = \floor{M(k,n)} + k \,,
 $$

\nin
  and
 $$
   \kaps(k,n) = \lams(k,n) = \dg(k,n) = \floor{M(k,n)} \,.
 $$
\etm

\bpf
In a packing of Type 1, all graphs are complete, and we have
 $\kaps(K_p) = \kappa(K_p) = p-1 = \Mad(K_p)$
 and $\omega(K_p) = p = \Mad(K_p) + 1$
 (and of course the values for $K_{p+1}$ are larger by eaxctly 1).
In case of Type 2, with $G\in \cG_{p,r}$ being the single non-complete
 graph in the packing (where $0<r<p$), we have $p-1<\Mad(G)<p$.
Then $\kaps(G) = \floor{\Mad(G)}$ and $\omega(G) = \floor{\Mad(G)} + 1$.
Summing up for $G_1,\dots,G_k$, the theorem follows.
\epf

As a matter of fact, the behavior of several functions is very similar.

\btm
For every parameter\/ $\pp\in\{\kaps,\lams,\dg\}$
 we have
 $$
   ( 1 - \epsilon_k ) \sqrt{k}\,n \leq \pp(k,n) \leq M(k,n) < \sqrt{k}\,n
 $$
 for every fixed\/ $k$, where\/ $\epsilon_k\to 0$ as\/ $k\to\infty$.
\etm

\bpf
The upper bound $\sqrt{k}$ is valid by Corollary \ref{c:M-sqrt-k}.
Since $\col(G)-1 = \dg(G)\geq \lams(G)\geq \kaps(G)$ holds
 for every graph $G$, it remains to prove the lower bound
 $\kaps(k,n) \geq ( 1 - \epsilon_k ) \sqrt{k}\,n$.

For this, we apply the construction by which we derived the
 asymptotic estimates of Theorem \ref{t:sqrt-k} along similar lines.
Let $k'\leq k$ be the largest integer of the form $k'=q^2+q+1$, where $q$
 is a prime power, and $n'$ be the largest multiple of $k'$ with $n'\leq n$.
Then $K_{n'}$ admits a packing $G_1,\dots,G_{k'}$ of $k'$ edge-disjoint
 copies of the complete $(q+1)$-partite graph
 $K_{\frac{n'}{k'},\frac{n'}{k'},\dots,\frac{n'}{k'}}$.
If fewer than $\frac{q}{q^2+q+1}\,n'$ vertices are removed from it, then
 more than one vertex class contains a non-deleted vertex.
Picking two of those vertices from distinct classes, they are adjacent and
 dominate the rest of $K_{\frac{n'}{k'},\frac{n'}{k'},\dots,\frac{n'}{k'}}$.
Thus, $\kaps(G_i) \geq \frac{q}{q^2+q+1}\,n'$ holds for all $q\leq i\leq k'$.

Distribute the edges of $E(K_n)\smin(\cup_{i=1}^{k'} E(G_i))$ among
 the current graphs $G_i$ and the further $k-k'$ graphs $G_{k'+1},\dots,G_k$
 to be constructed.
Using the facts $k-k'=o(k)$ and $n-n'<k$, we obtain

 \begin{eqnarray}
  \kaps(k,n) & \geq & \sum_{i=1}^{k'} \kaps(G_i) \ \geq \
   \frac{q}{q^2+q+1}\,k'n' \ = \ qn'  \nonumber \\
  & > & (\sqrt{k'}-1)\,n' \ \geq \ ( 1 - \epsilon_k ) \sqrt{k}\,n \,.  \nonumber
 \end{eqnarray}
\epf

From the first three parts of Theorem \ref{t:F-etal} one can see that
 the behavior of $\pp\in\{\omega,\chi,\ch\}$ is different because for
 those parameters we have $\pp(k,n)\leq n + f(k,n)$ with either
 $f(k,n)=c(k)$ (i.e., the error term depends on $k$ but is independent
 of $n$) or $f(k,n)=c(k)\cdot h(n)$ where $h(n)$ is sublinear.
On the other hand, our $\sqrt{k}\,n$ upper bound very quickly becomes
 substantially better than the corresponding estimates in the first three
 parts of Theorem \ref{t:F-etal}.

\btm   \label{t:omega-etc}
$\omega(k,n) \leq \chi(k,n) \leq \ch(k,n) < \sqrt{k}\,n + k$.
\etm

\section{Concluding remarks and open problems}
\label{s:concl}

In this section we present natural problems and directions for future
 research as a continuation suggested by the current research.

\bpm
Determine\/
 $\pp\ssl(k,N)$ for graph parameters\/ $\pp$ other than\/ $\pp(G)=\Mad(G)$.
\epm

\bpm
Study the product version of\/ $\pp(k,n)$ and\/ $\pp\ssl(k,N)$.
\epm

In cases where $\pp(k,n)$ is verified by a decomposition whose $\pp$-values
 form a balanced sequence of equal or nearly equal integers, Jensen's
 inequality is applicable and provides a strong upper bound.
This is the case with the function $\pp\ssl(k,N)$ as well:
 $$
   \max_{(G_1,\ldots,G_k) \ \mathrm{decompose} \ K_n}
     \Bigl\{ \, \prod_{i=1}^k \Mad(G_i) \Bigr\}
    \leq \Bigl( \frac{M(k,n)}{k} \Bigr)^k
    \leq \Bigl( \frac{M\ssl(k,\binom{n}{2})}{k} \Bigr)^k ,
 $$

 $$
   \max_{e(G_1)+\ldots+e(G_k) = N} \,
     \Bigl\{ \, \prod_{i=1}^k \Mad(G_i) \Bigr\}
    \leq \Bigl( \frac{M\ssl(k,N)}{k} \Bigr)^k
     < \bigl( 2N/k \bigr) ^{k/2} .
 $$
Our several results above show that this method works in certain
 ranges of $k,n$ for $M(k,n)$.
For instance, the following holds.

\btm
Let\/ $k = \binom{n}{2}/\binom{p}{2}$, and assume that a Steiner system\/
 $S(2,p,n)$ exists.
Then
 $$
   \max_{(G_1,\ldots,G_k) \ \mathrm{decompose} \ K_n}
     \Bigl\{ \, \prod_{i=1}^k \Mad(G_i) \Bigr\} =
    (p-1)^k = (p-1)^{ \frac{n(n-1)}{p(p-1)} } .
 $$
\etm

\bpf
On applying part $(i)$ of Theorem \ref{t:Stein} we see that in the
 present case $M(k,n)/k = p-1$.
The existence of Steiner systems ensures that the product of $\Mad$-values
 can now attain the universal upper bound $\bigl( M(k,n)/k \bigr)^k$.
\epf

In a similar way, constructions from multisets $\{a*K_p \cup b*K_{p+1}\}$
 yield exact ansewrs.
However, this kind of tight constructions requires that $k$ is relatively
 large with respect to $n$.
For $k$ small and $n$ large, the upper bound obtained via the
 geometric-arithmetic inequality is probably far from being tight,
 therefore the cases of fixed $k$ deserve attention.
In particular, it is natural to ask:

\bpm
Determine\/ \,$\max \{\Mad(G) \cdot \Mad(\ovG)\}$ \,over the graphs\/ $G$ of order\/ $n$.
In particular, is this quantity maximized by a pair\/
 $G = K_x$, \,$\ovG = K_n - K_x$, similarly to
 $\max \{\Mad(G) + \Mad(\ovG)\}$\,?
\epm

Concerning $\displaystyle{ \limsup_{n\to\infty} \max_{|G|=n}}\,
 \{\Mad(G) \cdot \Mad(\ovG)/n^2\}$,
Theorem \ref{t:k=2} implies the upper bound $25/64 = 0.390625$, and the
 construction in its proof yields $3/8 = 0.375$ as lower bound
 on $\liminf$.
We can improve on the latter, and obtain an estimate which is close to
 the $25/64$ bound within relative error less than $1.5$\,\%.

\bpn
$\displaystyle{ \liminf_{n\to\infty} \max_{|G|=n}}\, \{\Mad(G) \cdot \Mad(\ovG)/n^2\} \geq 2/\!\sqrt{27} > 0.3849$.
\epn

\bpf
Let $G\cong K_x$, where $x\approx n/\!\sqrt{3}$.
Then $\Mad(G)=x-1$, and $\Mad(\ovG) =
 \bigl(x\cdot(n-x) + (n-x)\cdot(n-1)\bigr)/n \approx
  (n-x)(n+x) / n= n\cdot \bigl(1-(x/n)^2\bigr) \approx 2n/3$.
\epf

\bpm
Improve on the lower bounds given in Proposition \ref{p:3--7} for\/
 $3\leq k\leq 7$, or prove that those estimates are optimal.
\epm

\bpm
Extend the formula derived in Theorem \ref{t:PSTS} from the range\/
 $k\geq \frac{1}{6} \, n(n-1)$ to\/ $k\geq \frac{1}{12} \, n(n-1)$.
\epm

Part of a solution would require the design of partial Steiner systems
 $PS(2,4,n)$ whose leave graphs admit edge decompositions into triangles.
An $S(2,4,n)$ system exists if and only if $n\equiv 1$ or 4 (mod 12).
However, a $\{K_3,K_4\}$-design requires that $a$ triangles and
 $b$ quadruples incident with any one vertex satisfy $2a+3b=n-1$, which
 puts considerable arithmetic limitations on the combination of numbers of
 $K_3$ and $K_4$ subgraphs occurring in a decomposition.
Therefore $M(k,n)$ in this range of $k$ will usually involve non-complete
 subgraphs in an optimal edge decomposition of $K_n$.
We note that for $n\equiv 1$ (mod 12) there exist both $S(2,3,n)$ and
 $S(2,4,n)$ systems, starting from $n=13$.
Also, $K_{10}$ decomposes into $\{9 * K_3 , \, 3 * K_4\}$, and
 $K_{12}$ decomposes into $\{4 * K_3 , \, 9 * K_4\}$.

\bpm
Characterize the extremal graphs for\/ $M(2,n)$ determined in Theorem \ref{t:k=2}.
\epm

Some results of Section \ref{ss:equiv} maybe useful in this direction also.

\bpm
For\/ $\pp\in\{\omega,\chi,\ch\}$ improve the\/ $\sqrt{k}\,n$ upper bound in
 Theorem~\ref{t:omega-etc}, and determine a strong lower bound on\/
 $\pp(k,n)$ as\/ $k=f(n)$ grows not very slowly.
\epm

In particular, concerning $\omega$, we gratefully
 thank a referee for advising us the following observation,
 which expresses equivalence between the Zaran\-kie\-wicz problem
 and the determination of $\omega(k, n)$.

\bpn
For all\/ $1\leq k\leq \binom{n}{2}$, $\omega(k, n)$ is equal to
 the largest number of edges in a\/ $C_4$-free bipartite graph
 with\/ $k$ and\/ $n$ vertices in its partite sets.
\epn

\bpf
Assume first that there exists a $C_4$-free bipartite graph
 $B\subset K_{k,n}$ with $m$ edges, with vertex classes
 $\{a_1,\dots,a_k\}$ and $\{b_1,\dots,b_n\}$.
Then the neighborhoods $N(a_1),\dots,N(a_k)$ form a system of
 $k$ sets, any two of which share at most one vertex.
Hence this can be viewed as an edge-disjoint packing of $k$
 complete subgraphs in the complete $K_n$ with vertex set
 $\{b_1,\dots,b_n\}$.
Of course, they can be extended to an edge decomposition
 $G_1,\dots,G_k$ of $K_n$.
Consequently, $\omega(k, n)\geq \sum_{i=1}^k
 \omega(G_i)\geq \sum_{i=1}^k |N(a_i)| = m$.

To derive the reverse inequality, assume that $G'_1,\dots,G'_k$
 is a decomposition of $K_n$ with $V(K_n)=\{b_1,\dots,b_n\}$,
  such that $m':= \sum_{i=1}^k \omega(G'_i) = \omega(k, n)$.
Choose a largest complete subgraph in $G'_i$ ($i=1,\dots,k$), say
 on the vertex subset $X_i\subset V(G'_i)$, $|X_i|=\omega(G'_i)$.
Since the $G'_i$ are edge-disjoint, any two of the $X_i$ have
 at most one vertex in common.
Let $B'\subset K_{k,n}$ be the bipartite graph whose vertex classes are
 $\{a_1,\dots,a_k\}$ and $\{b_1,\dots,b_n\}$, and each $a_i$ has
 neighborhood $N(a_i)=X_i$.
This $B'$ has $e(B') = \sum_{i=1}^k  |X_i|
  = \sum_{i=1}^k \omega(G_i) = m' = \omega(k, n)$,
 and it is $C_4$-free because each edge of $K_n$ belongs to
 only one $G_i$.
\epf

As we noted earlier, in the quadratic range of $k$ one can improve
 $\sqrt{k}\,n$ to  $(1-c)\sqrt{k}\,n$ with a suitably chosen constant
 $c>0$.
But the gap between $\pp(k,n)$ and $M(k,n)$ may be considerably larger
 than that.

Based on Propositions \ref{p:tree} and \ref{p:split}, and the subsequent
 related
 observations under them, the following question is very natural to ask.

\bpm
Characterize the graphs for which\/ $\Mad(G) = \ed(G)$ holds.
\epm

There are chances for a structural characterization because equality
 can be decided in polynomial time.
Namely, it is immediate to compute $\ed(G)$ in linear time, and although
 definitely not so easy, also $\Mad(G)$ can be determined by a
 polynomial-time algorithm as designed by Goldberg \cite{G-84}.
On the other hand, a characterization in terms of forbidden induced
 subgraphs does not exist.
This follows e.g.\ from the short discussion on the
 {\tt stackexchange.com} site \cite{Kor-22}.
Upon a question posed by Cyriac Antony there, Tuukka Korhonen (a.k.a.\
 Laakeri) commented that the complete join $H'=H+\overline{K}_{|H|}$ 
 satisfies $\Mad(H')=d(H')$ for every graph $H$.

\paragraph{Acknowledgement.}

We thank the referees for their very constructive remarks and suggestions.
Research of the second author was supported in part by the
 ERC Advanced Grant ``ERMiD''.


\newpage

\section*{Appendix: Alternative approach to Theorems \ref{t:list-k-n} and \ref{t:Mad-k}}
\addcontentsline{toc}{section}{Appendix: Alternative approach to Theorems \ref{t:list-k-n} and \ref{t:Mad-k}}

In this Appendix we provide an argument of algorithmic nature as an
 alternative proof of Theorem \ref{t:Mad-k}.
It demonstrates that every decomposition of $K_n$ can be transformed
 step by step to a unique well-defined representative list of graphs.
During this process the sum $\summad$ never decreases, and finally yields
 the claimed upper bound on $M(k,n)$.
The method is equally fine for a proof of Theorem \ref{t:list-k-n},
 hence determining $M\ssl(k,N)$ for arbitrary integers $N\geq k\geq 2$.

For convenience let us recall that the following will be proved:

\bsk

\nin
{\bf Theorem \ref{t:Mad-k}.}
{\it Let\/ $n\geq 3$ and\/ $2\leq k\leq \binom{n}{2}$ be given, and choose
  the integer\/ $p$ for which\/
  $k \binom{p}{2} \leq \binom{n}{2} < k \binom{p+1}{2}$ holds.
Let\/
 $q,r$ be the nonnegative uniquely determined integers such that\/
 $\binom{n}{2} = k \binom{p}{2} + qp + r$, where\/
 $0\leq q < k$ and $0\leq r < p$. Then

$$
  M(k,n) \leq \begin{cases}
   \begin{tabular}{ll}
    $kp - k + q$ & if\/  $r\leq (p-1)/2$\,; \\
    $kp - k + q + 1 - 2(p-r)/(p+1)$ & if\/  $r\geq (p-1)/2$\,.
   \end{tabular}
    \end{cases}
$$
Moreover, equality holds whenever\/ $K_n$ admits an edge decomposition
 into
  \begin{itemize}
   \item $q$ copies of\/ $K_{p+1}$ and\/ $k-q$ copies of\/ $K_p$, if $r=0$; or
   \item $q$ copies of\/ $K_{p+1}$, $k-q-1$ copies of\/ $K_p$,
    and one copy of any one graph\/ $G\in\fpr$, if\/ $0<r<p$.
  \end{itemize}
}

\bpf
Consider any edge decomposition $(G_1,\dots,G_k)$ of $K_n$.
We
 forget that the $G_i$ originate from a decomposition;
 we only view them as a list of individual graphs, and manipulate
 with the graphs themselves.
In this way, if we only assume that the starting collection of graphs
 is just a list $\sL=\{G_1,\dots,G_k\}$ with $e(G_1)+\ldots+e(G_k)=N$,
 then the process also leads to the determination of $M\ssl(k,N)$
 with the formula stated in Theorem \ref{t:list-k-n}.
However, in the sequel we restrict attention to $N=\binom{n}{2}$.

\msk

1.\quad
We replace each $G_i$ with $H_i$, where $H_i$ is the extremal
 representative of the graphs with $e(G_i)$ edges.
Due to Theorem \ref{t:max-g-m}, we have $\Mad(H_i)\geq \Mad(G_i)$.

\msk

2.\quad
We qualify the graphs $H_i$ into three categories, where the
 involved integer $p$ is uniquely determined by the number
 $e(H_i)$ of edges, and may vary for
 different subscripts among the $k$ graphs under consideration.
 \begin{description}
  \item[Type A:] $e(H_i)=\binom{p_i}{2}$.
  \item[Type B:] $\binom{p_i}{2} < e(H_i) < \binom{p_i+1}{2}$ and
   $\Mad(H_i) = 2e(H_i)/(p_i+1)$.
  \item[Type C:] $\binom{p_i}{2} < e(H_i) < \binom{p_i+1}{2}$ and
   $\Mad(H_i) = p_i-1$.
 \end{description}
At any stage of the procedure, we denote by $\cA$, $\cB$, $\cC$
 the actual collection (multiset) of graphs $H_i$ of
  type A, B, C, respectively.
If $e(H_i) = \frac{1}{2}(p_i+1)(p_i-1) = \binom{p_i}{2} + \frac{p_i-1}{2}$,
 then $H_i$ is qualified for both $\cB$ and $\cC$.
In this case we assume $H_i\in\cB$.

Beside these, we maintain an integer $s$ for the number of
 so-called ``spare edges'', which is initialized as $s=0$.

\msk

3.\quad
The following steps are performed, keeping $s+\sum_{i=1}^k e(H_i)$
 unchanged throughout, hence maintaining the sum equal to
  $\binom{n}{2}$.
 \begin{enumerate}[{3}.a]
  \item 
   First we eliminate the members of $\cC$.
  If a $H_i\in \cC$ occurs, say $H_i\cong \gpir$, we re-define
   $H_i:=K_{p_i}$ and $s:=s+r$.
    This keeps $\Mad(H_i)$ unchanged, and increases $s$.
   In this way all $H_i\in\cC$ are modified to type-A graphs
    at the beginning, and later on if a
    type-C graph arises, it is also eliminated immediately.
  \item \label{case:3b}
    Next we eliminate the members of $\cB$, possibly with one exception.
   If $\cB$ contains $H_i,H_j$ such that $H_i\cong\gpir$ and
    $H_j\cong \gpjr$, where either $p_i>p_j$, or $p_i=p_j$ and
    $r_i\leq r_j$, then transfer $r':=\min\{r_i,p_j-r_j\}$ edges
     from $H_i$ to $H_j$.
    We do the same if $|\cB|=1$, say $\cB=\{H_j\}$ where $H_j=\gpjr$,
     and $\cA$ has a member $H_i\cong K_{p_i+1}=G_{p_i,p_i}$
      for some $p_i>p_j$.
    
    Note that either $r'=p_j-r_j$ or $r'=p_j-r_j-1$ holds,
     because the assumptions $H_i,H_j\in\cB$ and $p_i\geq p_j$
     imply $p_j-r_j\leq (p_j+1)/2 \leq (p_i+1)/2 \leq r_i+1$.
   Hence, the following outcomes are possible after such a step.
    \begin{itemize}
     \item If $r'=r_i=p_j-r_j$, then $H_i\cong K_{p_i}$ and
      $H_j\cong K_{p_j+1}$, and both $H_i$ and $H_j$ have been
      moved from $\cB$ to $\cA$.
     \item If $r'=r_i=p_j-r_j-1$, then $H_i\cong K_{p_i}$ (moved from
      $\cB$ to $\cA$), and
      $H_j\cong G_{p_j,p_j-1}$ (hence remains a member of $\cB$).
     \item If $r'=p_j-r_j<r_i$, then $H_j\cong K_{p_j+1}$ (moved from
      $\cB$ to $\cA$), and $H_i\cong G_{p_i,r_i-r'}$, which either
      remains in $\cB$ or has been moved from $\cB$ to $\cC$.
       In the latter case it is then immediately transformed to
      $K_{p_i}$ by step 3.a, hence becomes a member of $\cA$.
 \end{itemize}

   In either of these three cases, the value of $\Mad(H_j)$
    has been increased by exactly $2r'\!/(p_j+1)$.
   On the other hand, the decrease of $\Mad(H_i)$ is at most
    $2r'\!/(p_i+1)$.
   Indeed, removing the $r'$ edges from $H_i$ one by one, until
    reaching the threshold of $r_i=(p_i-1)/2$ the omission of each edge
    yields a loss of $2/(p_i+1)$, and no decrease occurs afterwards.
   Thus, the overall increase of $\Mad(H_i)+\Mad(H_j)$ is
    not smaller than $\frac{2(p_i-p_j)r'}{(p_i+1)(p_j+1)}$.
    
   In case $p_i=p_j$ it is possible that $\Mad(H_i)+\Mad(H_j)$
    remains unchanged.
   This happens if the modified $H_j$ is $K_{p_i+1}$ while $H_i$
    still remains a member of $\cB$
    (what means $r_i+r_j-p_i\geq (p_i-1)/2$).
   
   Performing 3.a and 3.b as long as possible, we obtain $\cC=\es$
    and either $\cB=\es$ or $|\cB|=1$.

  \item The next goal is to make a balance in $\cA \cup \cB$.
   Recall that $G_{p',p'}\cong K_{p'+1}$ holds for every $p'\geq 1$.
   Assume that $H_i=G_{p_i,p_i}$ is the largest member of $\cA$, and
    $H_j=G_{p_j,p_j}$ is its smallest member; moreover if $\cB\neq \es$,
    let $\cB=\{\gpkr\}$ (where $0<r_k<p_k$).
   If $\max\{p_i,p_j,p_k\} > \min\{p_i,p_j,p_k\}+1$ (disregarding $p_k$
    if $\cB=\es$), select the two of
    $H_i,H_j,H_k$ whose $p$-value is largest and smallest, respectively,
    and perform the transformation 3.\ref{case:3b} on them.
   (If the step is applied on $H_i,H_j\in \cA$, then $s$ is updated
    exactly as $s:=s+p_i-p_j-1$.)
   Repeat this until all $p$-values in $\cA \cup \cB$ differ by at most 1.
   In the same way as in the analysis of 3.b, it is readily seen
    that the sum of $\Mad$-values has been increased (and also $s$).
   
  \item In the last phase the $s$ spare edges are eliminated.
   If $\cB\ne\es$ and $H_k\cong\gpkr$, first $r'':=\min\{s,p_k-r_k\}$
    edges are inserted in $H_k$.
   The process is finished if $r''=s$.
    Otherwise we update $s:=s-r''$; at that moment $\cB=\cC=\es$, and the
     family $\cA$ as a multiset is of the form
      $\{a * K_{p'} , b * K_{p'+1}\}$ (meaning $a$ copies of $K_{p'}$
      and $b$ copies of $K_{p'+1}$)
       for some $p'\geq 2$, $a\geq 1$, $b\geq 0$.
  As long as $s\geq p'$, pick a copy of $K_{p'}$ from $\cA$,
   enlarge it to $K_{p'+1}$, and update $s:=s-p'$.
    Once this has been done for all $H_i\cong K_{p'}$, and still
     $s>0$ holds, update $p':=p'+1$ and continue in the same way.
   
  \item If $s=0$, stop. If $s>0$, replace one copy of $K_{p'}$ with
   the graph $G_{p',s}$.
    Call the current family $H_1,\dots,H_k$ the representative of
   the pair $(k,n)$.
   
   Due to the choice of $p,q,r$ with respect to $k$,
    at this point $\cA\cup \cB\cup \cC$ is of the form either
      $\{a * K_{p} , b * K_{p+1} , 1 * G_{p,r}\}$
       with $1\leq r<p$ and $a,b\geq 0$ ($a=0$ or $b=0$ is possible,
       $a+b=k-1$),
        or
      $\{a * K_{p} , b * K_{p+1}\}$
        (here $a+b=k$, possibly with $b=0$), where $a,b$ depend on
        $n$ and $k$, and are easily expressible in terms of $q$ and $r$,
      as stated in the theorem.
 \end{enumerate}

Observe that, starting from an arbitrary edge decomposition of $K_n$
 into $k$ subgraphs $G_1,\dots,G_k$, the procedure terminates
 with the representative of the pair $(k,n)$ as the unique final
  family $H_1,\dots,H_k$ determined above.
Thus,

 $$
   M(k,n) \leq \Mad(H_1) + \ldots + \Mad(H_k) \,.
 $$
The numerical formula follows from the facts that $\Mad(K_p)=p-1$
 (applied also as $\Mad(K_{p+1})=p$) and that $\Mad(G_{p,r})$ is
 either $p$ or $(p^2-p+2r)/(p+1)$, depending on
 whether $r\leq (p-1)/2$ or bigger.
Namely, the procedure ends with the combination of graphs as
 displayed in the following table:

\msk

\begin{center}
\begin{tabular}{lllll}
\quad graphs $H_i$ && total of $\Mad(H_i)$ && side condition \\ \hline
$q$ copies of $K_{p+1}$ && $qp$ && --- \\
$k-q-1$ copies of $K_p$ && $(k-q-1)(p-1)$ && --- \\
$1$ copy of $G_{p,r}$ && $p-1$ && $r\leq (p-1)/2$ \\
$1$ copy of $G_{p,r}$ && $p - 2(p-r)/(p+1)$ && $r > (p-1)/2$
\end{tabular}
\end{center}

\msk

\nin
Thus, summing up the entries in the middle column according to the
 side condition depending on the value of $r$, the
 claimed upper bound follows.

The cases of equality, as stated in the last part of the theorem,
 can also be verified.
Indeed, if an assumed edge decomposition of $K_n$ exists, then all
 computations above hold with equality, and the proved upper bound
 on $M(k,n)$ is attained.
\epf

\end{document}